# SPECTRAL MEASURE OF LARGE RANDOM HANKEL, MARKOV AND TOEPLITZ MATRICES[1]

By Włodzimierz Bryc, Amir Dembo and Tiefeng Jiang

*University of Cincinnati, Stanford University and University of Minnesota*

We study the limiting spectral measure of large symmetric random matrices of linear algebraic structure.

For Hankel and Toeplitz matrices generated by i.i.d. random variables $\{X_k\}$ of unit variance, and for symmetric Markov matrices generated by i.i.d. random variables $\{X_{ij}\}_{j>i}$ of zero mean and unit variance, scaling the eigenvalues by $\sqrt{n}$ we prove the almost sure, weak convergence of the spectral measures to universal, nonrandom, symmetric distributions $\gamma_H$, $\gamma_M$ and $\gamma_T$ of unbounded support. The moments of $\gamma_H$ and $\gamma_T$ are the sum of volumes of solids related to Eulerian numbers, whereas $\gamma_M$ has a bounded smooth density given by the free convolution of the semicircle and normal densities.

For symmetric Markov matrices generated by i.i.d. random variables $\{X_{ij}\}_{j>i}$ of mean $m$ and finite variance, scaling the eigenvalues by $n$ we prove the almost sure, weak convergence of the spectral measures to the atomic measure at $-m$. If $m=0$, and the fourth moment is finite, we prove that the spectral norm of $\mathbf{M}_n$ scaled by $\sqrt{2n\log n}$ converges almost surely to 1.

**1. Introduction and main results.** For a symmetric $n\times n$ matrix $\mathbf{A}$, let $\lambda_j(\mathbf{A})$, $1\le j\le n$, denote the eigenvalues of the matrix $\mathbf{A}$, written in a nonincreasing order. The spectral measure of $\mathbf{A}$, denoted $\hat{\mu}(\mathbf{A})$, is the empirical distribution of its eigenvalues, namely

$$\hat{\mu}(\mathbf{A}) = \frac{1}{n}\sum_{j=1}^{n}\delta_{\lambda_j(\mathbf{A})}$$

Received May 2004; revised January 2005.

[1]Supported in part by NSF Grants INT-0332062, DMS-00-72331, DMS-03-08151, DMS-04-49365 and DMS-05-04198.

*AMS 2000 subject classifications.* Primary 15A52; secondary 60F99, 62H10, 60F10.

*Key words and phrases.* Random matrix theory, spectral measure, free convolution, Eulerian numbers.

[so when $\mathbf{A}$ is a random matrix, $\hat{\mu}(\mathbf{A})$ is a random measure on $(\mathbb{R}, \mathcal{B})$].

Large-dimensional random matrices are of much interest in statistics, where they play a pivotal role in multivariate analysis. In his seminal paper, Wigner [24] proved that the spectral measure of a wide class of symmetric random matrices of dimension $n$ converges, as $n \to \infty$, to the semicircle law (also called the Sato–Tate measure, see [21] and the references therein). Much work has since been done on related random matrix ensembles, either composed of (nearly) independent entries, or drawn according to weighted Haar measures on classical (e.g., orthogonal, unitary, simplectic) groups. The limiting behavior of the spectrum of such matrices and their compositions is of considerable interest for mathematical physics (see [17] and the references therein). In addition, such random matrices play an important role in operator algebra studies initiated by Voiculescu, known now as the free (noncommutative) probability theory (see [12] and the many references therein). The study of large random matrices is also related to interesting questions of combinatorics, geometry and algebra (see [9], or, e.g., [22]). In his recent review paper [1], Bai proposes the study of large random matrix ensembles with certain additional linear structure. In particular, the properties of the spectral measures of random Hankel, Markov and Toeplitz matrices with independent entries are listed among the unsolved random matrix problems posed in [1], Section 6. We shall provide here the solution for these three problems.

We note in passing that Hankel matrices arise, for example, in polynomial regression, as the covariance for the least squares parameter estimation for the model $\sum_{i=0}^{p-1} b_i x^i$, observed at $x = x_1, \ldots, x_n$ in the presence of additive noise (see [20], page 36). Toeplitz matrices appear as the covariance of stationary processes, in shift-invariant linear filtering, and in many aspects of combinatorics, time series and harmonic analysis. See [10] for classical results on deterministic Toeplitz matrices, or [7] and the references therein, for their applications to certain random matrices. The infinitesimal generators of continuous-time Markov processes on finite state spaces are given by matrices with row-sums zero (which we call Markov matrices). Such matrices also play an important role in graph theory, as the Laplacian matrix of each graph is of this form, with its eigenvalues related to numerous graph invariants; see [15].

We next specify the corresponding ensembles of random matrices studied here. Let $\{X_k : k = 0, 1, 2, \ldots\}$ be a sequence of i.i.d. real-valued random vari-ables. For $n \in \mathbb{N}$, define a random $n \times n$ Hankel matrix $\mathbf{H}_n = [X_{i+j-1}]_{1 \le i,j \le n}$,

$$\mathbf{H}_n = \begin{bmatrix} X_1 & X_2 & \cdots & \cdots & X_{n-1} & X_n \\ X_2 & X_3 & & & X_n & X_{n+1} \\ \vdots & & & & X_{n+1} & X_{n+2} \\ & & & \iddots & & \\ X_{n-2} & X_{n-1} & & & & \vdots \\ & & \iddots & & & \\ X_{n-1} & X_n & & & X_{2n-3} & X_{2n-2} \\ X_n & X_{n+1} & \cdots & \cdots & X_{2n-2} & X_{2n-1} \end{bmatrix}, \tag{1.1}$$

and a random $n \times n$ Toeplitz matrix $\mathbf{T}_n = [X_{|i-j|}]_{1 \le i,j \le n}$,

$$\mathbf{T}_n = \begin{bmatrix} X_0 & X_1 & X_2 & \cdots & X_{n-2} & X_{n-1} \\ X_1 & X_0 & X_1 & & & X_{n-2} \\ X_2 & X_1 & X_0 & & \ddots & \vdots \\ \vdots & & & \ddots & & X_2 \\ X_{n-2} & & & & X_0 & X_1 \\ X_{n-1} & X_{n-2} & \cdots & X_2 & X_1 & X_0 \end{bmatrix}. \tag{1.2}$$

The limiting spectral distribution for a Toeplitz matrix $\mathbf{T}_n$ is as follows.

THEOREM 1.1. *Let $\{X_k : k = 0, 1, 2, \ldots\}$ be a sequence of i.i.d. real-valued random variables with* $\mathrm{Var}(X_1) = 1$. *Then with probability 1, $\hat{\mu}(\mathbf{T}_n/\sqrt{n})$ converges weakly as $n \to \infty$ to a nonrandom symmetric probability measure $\gamma_T$ which does not depend on the distribution of $X_1$, and has unbounded support.*

The spectrum of nonrandom Toeplitz matrices, the rows of which are typically absolutely summable, is well approximated by its counterpart for circulant matrices (cf. [10], page 84). In contrast, note that the limiting distribution $\gamma_T$ is not normal as the calculation shows that the fourth moment is $m_4 = 8/3$. This differs from the analogous results for random circulant matrices (see [4]), a fact that has been independently noticed also in references [3, 11].

Our next result gives the limiting spectral distribution for a Hankel matrix $\mathbf{H}_n$.

THEOREM 1.2. *Let $\{X_k : k = 0, 1, 2, \ldots\}$ be a sequence of i.i.d. real-valued random variables with* $\mathrm{Var}(X_1) = 1$. *Then with probability 1, $\hat{\mu}(\mathbf{H}_n/\sqrt{n})$*

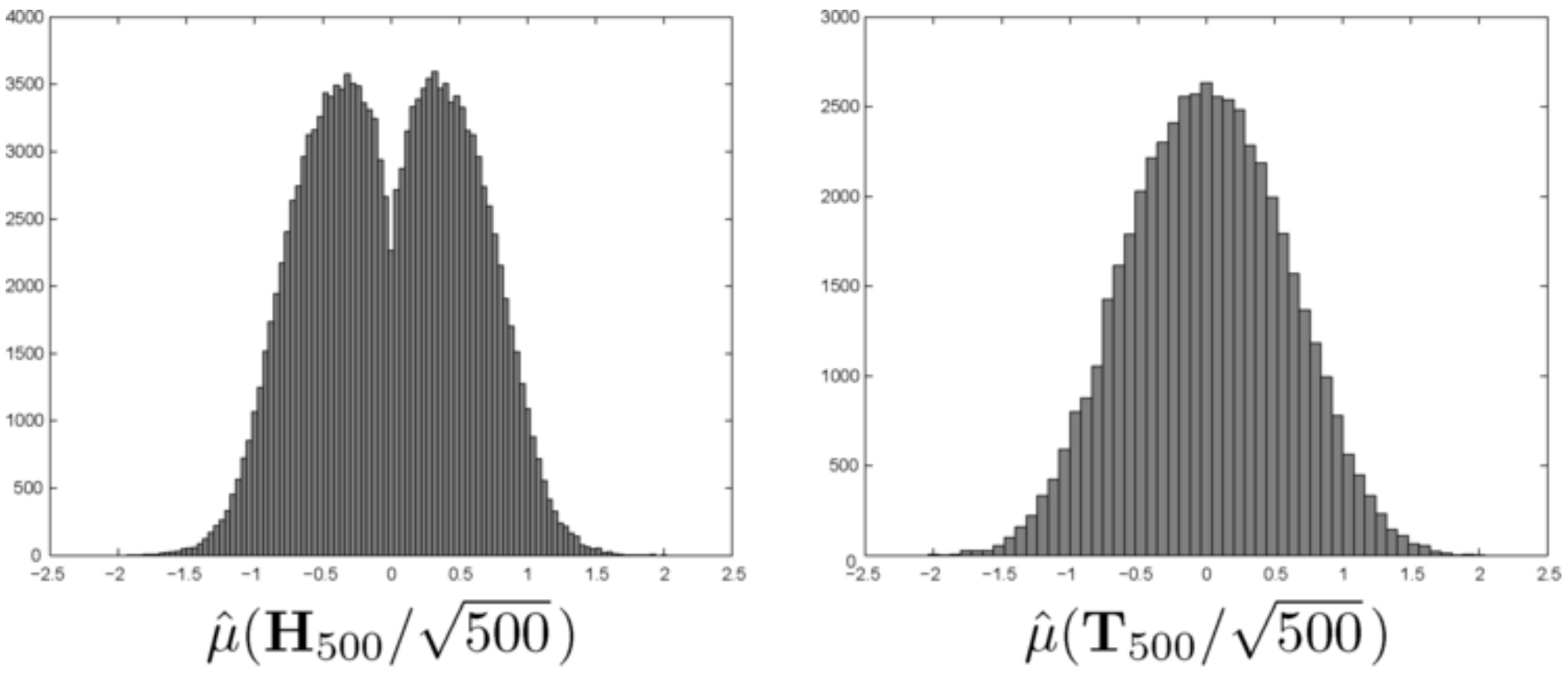

FIG. 1. *Histograms of the empirical distribution of eigenvalues of* 100 *realizations of the Hankel and Toeplitz matrices with standardized triangular* $U - U'$ *entries.*

*converges weakly as* $n \to \infty$ *to a nonrandom symmetric probability measure* $\gamma_H$ *which does not depend on the distribution of* $X_1$, *has unbounded support and is not unimodal.*

(Recall that a symmetric distribution $\nu$ is said to be unimodal, if the function $x \mapsto \nu((-\infty, x])$ is convex for $x < 0$.)

REMARK 1.1. Theorems 1.1 and 1.2 fall short of establishing that the limiting distributions have smooth densities and that the density of $\gamma_H$ is bimodal. Simulations suggest that these properties are likely to be true; see Figure 1.

REMARK 1.2. Consider the empirical distribution of singular values of the nonsymmetric random $n \times n$ Toeplitz matrix $\mathbf{R}_n = [X_{i-j}]_{1 \le i,j \le n}$. It follows from Theorem 1.2 that as $n \to \infty$, with probability 1, $\hat{\mu}((\mathbf{R}_n \mathbf{R}_n^T)^{1/2}/\sqrt{n}) \to \nu$ weakly, where $\nu([0,x]) = \gamma_H([-x,x]), x > 0$. Indeed, let $\mathbf{J}_n = [\mathbb{1}_{i+j=n+1}]_{1 \le i,j \le n}$, noting that $\mathbf{J}_n \times \mathbf{R}_n^T$ is the Hankel matrix $\mathbf{H}_n$ for $\{X_{k-n} : k = 0, 1, \ldots\}$ to which Theorem 1.2 applies. Since $\mathbf{J}_n^2 = \mathbf{I}_n$, and both $\mathbf{J}_n$ and $\mathbf{J}_n \times \mathbf{R}_n^T$ are symmetric, we have $\mathbf{R}_n \mathbf{R}_n^T = (\mathbf{R}_n \mathbf{J}_n)^T \mathbf{J}_n \mathbf{R}_n^T = \mathbf{H}_n^2$. Thus the singular values of matrix $\mathbf{R}_n$ are the absolute values of the (real) eigenvalues of the symmetric Hankel matrix $\mathbf{H}_n$.

We now turn to the Markov matrices $\mathbf{M}_n$. Let $\{X_{ij} : j \ge i \ge 1\}$ be an infinite upper triangular array of i.i.d. random variables and define $X_{ji} = X_{ij}$ for $j > i \ge 1$. Let $\mathbf{M}_n$ be a random $n \times n$ symmetric matrix given by

$$\mathbf{M}_n = \mathbf{X}_n - \mathbf{D}_n, \tag{1.3}$$

where $\mathbf{X}_n = [X_{ij}]_{1\le i,j\le n}$ and $\mathbf{D}_n = \operatorname{diag}(\sum_{j=1}^n X_{ij})_{1\le i\le n}$ is a diagonal matrix, so each of the rows of $\mathbf{M}_n$ has a zero sum (note that the values of $X_{ii}$ are irrelevant for $\mathbf{M}_n$), that is,

$$\mathbf{M}_n = \begin{bmatrix} -\sum_{j=2}^n X_{1j} & X_{12} & X_{13} & & \cdots & X_{1n} \\ X_{21} & -\sum_{j\ne 2}^n X_{2j} & X_{23} & & \cdots & X_{2n} \\ \vdots & & \ddots & & & \vdots \\ X_{k1} & X_{k2} & \cdots & -\sum_{j\ne k}^n X_{kj} & \cdots & X_{kn} \\ \vdots & \vdots & & & \ddots & \vdots \\ X_{n1} & X_{n2} & \cdots & & & -\sum_{j=1}^{n-1} X_{nj} \end{bmatrix}.$$

Wigner's classical result says that $\hat\mu(\mathbf{X}_n/\sqrt{n})$ converges weakly as $n\to\infty$ to the (standard) semicircle law with the density $\sqrt{4-x^2}/(2\pi)$ on $(-2,2)$. For normal $\mathbf{X}_n$ and normal i.i.d. diagonal $\tilde{\mathbf{D}}_n$ independent of $\mathbf{X}_n$, the weak limit of $\hat\mu((\mathbf{X}_n-\tilde{\mathbf{D}}_n)/\sqrt{n})$ is the free convolution of the semicircle and standard normal measures; see [17] and the references therein (see also [2] for the definition and properties of the free convolution). This predicted result holds also for the Markov matrix $\mathbf{M}_n$, but the problem is nontrivial because $\mathbf{D}_n$ strongly depends on $\mathbf{X}_n$.

THEOREM 1.3. *Let $\{X_{ij}: j\ge i\ge 1\}$ be a collection of i.i.d. random variables with $\mathbb{E}X_{12}=0$ and $\operatorname{Var}(X_{12})=1$. With probability 1, $\hat\mu(\mathbf{M}_n/\sqrt{n})$ converges weakly as $n\to\infty$ to the free convolution $\gamma_M$ of the semicircle and standard normal measures. This measure $\gamma_M$ is a nonrandom symmetric probability measure with smooth bounded density, does not depend on the distribution of $X_{12}$ and has unbounded support.*

If the mean of $X_{ij}$ is not zero, the following result is relevant.

THEOREM 1.4. *Let $\{X_{ij}: i,j\in\mathbb{N},\ j\ge i\ge 1\}$ be a collection of i.i.d. random variables with $\mathbb{E}X_{12}=m$ and $\mathbb{E}X_{12}^2<\infty$. Then $\hat\mu(\mathbf{M}_n/n)$ converges weakly to $\delta_{-m}$ as $n\to\infty$.*

Turning to the asymptotic of the spectral norm $\|\mathbf{M}_n\| := \max\{\lambda_1(\mathbf{M}_n), -\lambda_n(\mathbf{M}_n)\}$ of the symmetric matrix $\mathbf{M}_n$, that is, the largest absolute value of its eigenvalues, we have the following

THEOREM 1.5. *Let $\{X_{ij}: i, j \in \mathbb{N}, j \geq i \geq 1\}$ be a collection of i.i.d. random variables with $\mathbb{E}X_{12} = 0$, $\mathrm{Var}(X_{12}) = 1$ and $\mathbb{E}X_{12}^4 < \infty$. Then*

$$\lim_{n\to\infty} \frac{\|\mathbf{M}_n\|}{\sqrt{2n\log n}} = 1 \qquad a.s.$$

If the mean of $X_{ij}$ is not zero, the following result is relevant.

COROLLARY 1.6. *Suppose $\mathbb{E}X_{12} = m$ and $\mathbb{E}X_{12}^4 < \infty$. Then*

$$\lim_{n\to\infty} \frac{\|\mathbf{M}_n\|}{n} = |m| \qquad a.s.$$

Theorem 1.5 reveals a scaling in $n$ that differs from that of the spectral norm of Wigner's ensemble, where under the same conditions, almost surely,

$$\lim_{n\to\infty} \frac{\|\mathbf{X}_n\|}{\sqrt{n}} = 2 \tag{1.4}$$

(cf. [1], Theorem 2.12). As shown in Section 2 enroute to proving Theorems 1.4, 1.5 and Corollary 1.6, this is due to the domination of the diagonal terms of $\mathbf{M}_n$ in determining its spectral norm.

REMARK 1.3. The asymptotic of the spectral norm of random Toeplitz $\mathbf{T}_n$ and Hankel $\mathbf{H}_n$ matrices is not addressed in this work.

Theorems 1.4, 1.5 and Corollary 1.6 are proved in Section 2. The proofs of Theorems 1.1 and 1.2, which are similar to each other, ultimately rely on the method of moments and the well-known relation

$$\int x^k \hat{\mu}(\mathbf{A})(dx) = \frac{1}{n} \operatorname{tr} \mathbf{A}^k$$

for an $n \times n$ symmetric matrix $\mathbf{A}$. We begin in Section 3 by introducing the combinatorial structures which describe the moments of the limiting distributions. (Proofs of the properties of the limiting distributions are postponed to the Appendix.) Then in Section 4.1 we use truncation arguments to reduce the theorems to the case when the expected values of the moments of the spectral measures are finite. In Section 4.2 we show that under suitable integrability assumptions the expected values of moments of the spectral measures converge to the corresponding expressions from Section 3 as the size of the matrix $n \to \infty$. Representing the moments as traces, we use independence of the entries and combinatorial arguments to discard the irrelevant terms in the expansions (4.7) and (4.12). In Section 4.4 we show that the moments of the spectral measures are concentrated around their means, which allows us to conclude the proofs in Section 4.5.

The proof of Theorem 1.3 follows a similar plan, with truncation argument in Section 4.1, followed by combinatorial analysis of expansion (4.17) for the traces and concentration of moments in Section 4.6.

**2. Proofs of Theorems 1.4, 1.5 and Corollary 1.6.** We need the following result, which follows by Chebyshev's inequality from [18], Section 6, Theorem 5 or [19], Section 5, Corollary 5.

LEMMA 2.1 (Sakhanenko). *Let $\{\xi_i;\ i = 1, 2, \dots\}$ be a sequence of independent random variables with mean zero and $\mathbb{E}\xi_i^2 = \sigma_i^2$. If $\mathbb{E}|\xi_i|^p < \infty$ for some $p > 2$, then there exists a constant $C > 0$ and $\{\eta_i,\ i = 1, 2, \dots\}$, a sequence of independent normally distributed random variables with $\eta_i \sim N(0, \sigma_i^2)$ such that*

$$\mathbb{P}\left(\max_{1\le k\le n} |S_k - T_k| > x\right) \le \frac{C}{1 + |x|^p} \sum_{i=1}^{n} \mathbb{E}|\xi_i|^p$$

*for any $n$ and $x > 0$, where $S_k = \sum_{i=1}^{k} \xi_i$ and $T_k = \sum_{i=1}^{k} \eta_i$.*

PROOF OF THEOREM 1.5. Hereafter let $b(n) = \sqrt{2n \log n}$ denote the normalization function for Theorem 1.5.

It follows from (1.3) that $|\|\mathbf{M}_n\| - \|\mathbf{D}_n\|| \le \|\mathbf{X}_n\|$. So, by (1.4) and the definition of $\mathbf{D}_n$, it suffices to show that as $n \to \infty$,

$$W_n := \frac{1}{b(n)} \max_{i=1}^{n} \left\{ \left| \sum_{j=1}^{n} X_{ij} \right| \right\} \to 1 \qquad \text{a.s.} \tag{2.1}$$

We first show the upper bound, that is,

$$\limsup_{n\to\infty} W_n \le 1 \qquad \text{a.s.} \tag{2.2}$$

Note that $\{X_{ij};\ j \ge 1\}$ is a sequence of i.i.d. random variables for each $i \ge 1$. By Lemma 2.1 and the condition that $\mathbb{E}|X_{12}|^4 < \infty$, for each $i \ge 1$, there exists a sequence of independent standard normals $\{Y_{ij}; j \ge 1\}$ such that

$$\max_{i=1}^{n} \mathbb{P}\left( \max_{k=1}^{n} \left| \sum_{j=1}^{k} (X_{ij} - Y_{ij}) \right| > x \right) \le \frac{Cn}{x^4} \tag{2.3}$$

for all $x > 0$ and $n \ge 1$, where $C$ is a constant which does not depend on $n$ and $x$ (note that two sequences $\{Y_{ij}; j \ge 1\}$ for different values of $i$ are not independent of each other). We claim that

$$U_n := \frac{1}{b(n)} \max_{i=1}^{n} \left\{ \left| \sum_{j=1}^{n} (X_{ij} - Y_{ij}) \right| \right\} \to 0 \qquad \text{a.s.} \tag{2.4}$$

as $n \to \infty$. First,

$$\max_{k=2^m}^{2^{m+1}} U_k \le \frac{1}{b(2^m)} \max_{i=1}^{2^{m+1}} \max_{k=1}^{2^{m+1}} \left\{ \left| \sum_{j=1}^{k} (X_{ij} - Y_{ij}) \right| \right\}.$$

By (2.3), for any $\varepsilon > 0$,

$$\mathbb{P}\left(\max_{k=2^m}^{2^{m+1}} U_k \geq \varepsilon\right) \leq 2^{m+1}\mathbb{P}\left(\max_{k=1}^{2^{m+1}}\left|\sum_{j=1}^{k}(X_{ij} - Y_{ij})\right| \geq \varepsilon b(2^m)\right)$$
$$\leq \frac{C_\varepsilon}{m^2}$$

for some constant $C_\varepsilon$ depending only on $\varepsilon$. Since $\varepsilon > 0$ is arbitrary, by the Borel–Cantelli lemma, $\max_{k=2^m}^{2^{m+1}} U_k \to 0$ a.s. as $m \to \infty$, which implies (2.4). Let

$$V_n = \frac{1}{b(n)} \max_{i=1}^{n}\left|\sum_{j=1}^{n} Y_{ij}\right|.$$

By the definitions in (2.1) and (2.4), we have that $W_n \leq U_n + V_n$, so by (2.4) we get (2.2) as soon as we show that $\limsup_{n\to\infty} V_n \leq 1$. To this end, fix $\delta > 0$ and $\alpha > 1/\delta$. Then,

$$\begin{aligned} &\mathbb{P}\left(\max_{n=m^\alpha}^{(m+1)^\alpha} V_n \geq 1 + \delta\right) \\ &\quad \leq (m+1)^\alpha \mathbb{P}\left(\max_{n=1}^{(m+1)^\alpha}\left|\sum_{j=1}^{n} Y_{1j}\right| \geq (1+\delta)b(m^\alpha)\right) \\ &\quad \leq 2(m+1)^\alpha \mathbb{P}\left(\left|\sum_{j=1}^{(m+1)^\alpha} Y_{1j}\right| \geq (1+\delta)b(m^\alpha)\right), \end{aligned} \tag{2.5}$$

where Lévy's inequality is used in the second step. Since $Y_{ij}$'s are independent standard normals, $\xi := (m+1)^{-\alpha/2}\sum_{j=1}^{(m+1)^\alpha} Y_{1j}$ is a standard normal random variable. Thus, by the well-known normal tail estimate

$$\frac{1}{\sqrt{2\pi}}\frac{x}{1+x^2}e^{-x^2/2} \leq \mathbb{P}(\xi > x) \leq \frac{1}{\sqrt{2\pi}}\frac{1}{x}e^{-x^2/2} \qquad \text{for } x > 0, \tag{2.6}$$

we see that

$$\mathbb{P}(|\xi| \geq (1+\delta)(m+1)^{-\alpha/2}b(m^\alpha)) \leq \widehat{C}_\delta m^{-\alpha(1+\delta)}$$

for some constant $\widehat{C}_\delta > 0$. Consequently, for some $C'_\delta > 0$ and all $m$, by (2.5),

$$\mathbb{P}\left(\max_{n=m^\alpha}^{(m+1)^\alpha} V_n \geq 1 + \delta\right) \leq C'_\delta m^{-\alpha\delta}.$$

With $\alpha\delta > 1$, we have by the Borel–Cantelli lemma that

$$\limsup_{m\to\infty}\left\{\max_{n=m^\alpha}^{(m+1)^\alpha} V_n\right\} \leq 1 + \delta \qquad \text{a.s.}$$

It follows that $\limsup_{n\to\infty} V_n \le 1+\delta$ a.s. and taking $\delta \downarrow 0$ we obtain (2.2).

We next prove that

$$\liminf_{n\to\infty} W_n \ge 1 \qquad \text{a.s.} \tag{2.7}$$

To this end, fixing $1/3 > \varepsilon > \delta > 0$, let $n_\varepsilon := [n^{1-\varepsilon}] + 1$. Then,

$$\begin{aligned} W_n &\ge \frac{1}{b(n)} \max_{i=1}^{n_\varepsilon} \left| \sum_{j=1}^{n} X_{ij} \right| \\ &\ge \frac{1}{b(n)} \max_{i=1}^{n_\varepsilon} \left| \sum_{j=n_\varepsilon+1}^{n} X_{ij} \right| - \frac{1}{b(n)} \max_{i=1}^{n_\varepsilon} \left| \sum_{j=1}^{n_\varepsilon} X_{ij} \right| \qquad (2.8) \\ &=: V_{n,1} - V_{n,2}. \end{aligned}$$

By (2.2), $\limsup_{n\to\infty} W_{n_\varepsilon} \le 1$ a.s. Thus, with $b(n_\varepsilon)/b(n) \to 0$ as $n \to \infty$, we have that

$$V_{n,2} = W_{n_\varepsilon} \frac{b(n_\varepsilon)}{b(n)} \to 0 \qquad \text{a.s.} \tag{2.9}$$

Since $\{X_{ij};\ 1 \le i \le n_\varepsilon,\ n_\varepsilon < j \le n\}$ are i.i.d. for any $n \ge 1$, it follows that

$$\mathbb{P}(V_{n,1} \le 1 - 3\delta) = \mathbb{P}\left( \left| \sum_{j=1}^{n-n_\varepsilon} X_{1j} \right| \le (1-3\delta) b(n) \right)^{n_\varepsilon}. \tag{2.10}$$

With $b(n) \ge \sqrt{n}$, by Lemma 2.1 there exists a sequence of independent standard normals $\{Y_j\}$ such that for some $C = C(\delta) < \infty$ and all $n$

$$\mathbb{P}\left( \left| \sum_{j=1}^{n-n_\varepsilon} X_{1j} - \sum_{j=1}^{n-n_\varepsilon} Y_j \right| \ge \delta b(n) \right) \le C n^{-1}. \tag{2.11}$$

Further, by the left inequality of (2.6) we have that for all $n$ sufficiently large,

$$\begin{aligned} \mathbb{P}\left( \left| \sum_{j=1}^{n-n_\varepsilon} Y_j \right| \le (1-2\delta) b(n) \right) &\le \mathbb{P}(|Y_1| \le (1-\delta)\sqrt{2\log n}\,) \\ &\le 1 - 2n^{-(1-\delta)}. \end{aligned}$$

Combining this bound with (2.11) and (2.10) we get that for all $n$ large enough

$$\begin{aligned} \mathbb{P}(V_{n,1} \le 1 - 3\delta) &\le (1 - 2n^{-(1-\delta)} + C n^{-1})^{n_\varepsilon} \\ &\le (1 - n^{-(1-\delta)})^{n^{1-\varepsilon}} \le e^{-n^{\varepsilon-\delta}}. \end{aligned}$$

Recall that $\varepsilon > \delta$, implying that $\sum_{n\ge 1}\mathbb{P}(V_{n,1}\le 1-3\delta)<\infty$. By the Borel–Cantelli lemma,

$$\liminf_{n\to\infty} V_{n,1}\ge 1-3\delta \qquad \text{a.s.}$$

This together with (2.8) and (2.9) implies that almost surely $\liminf_{n\to\infty} W_n \ge 1-3\delta$, and the lower bound (2.7) follows by taking $\delta\downarrow 0$. □

Proof of Corollary 1.6. Let $\widetilde{\mathbf{M}}_n$ denote the Markov matrix obtained when $\widetilde{X}_{ij}=X_{ij}-\mathbb{E}X_{ij}$ replaces $X_{ij}$ in (1.3). Obviously,

$$\mathbf{M}_n=\widetilde{\mathbf{M}}_n+\mathbf{Y}_n, \tag{2.12}$$

where $\mathbf{Y}_n=[Y_{ij}]$ is the $n\times n$ matrix with $Y_{ij}=m-nm\mathbf{1}_{i=j}$. Clearly, $\lambda_1(\mathbf{Y}_n)=0$, $\lambda_2(\mathbf{Y}_n)=\cdots=\lambda_n(\mathbf{Y}_n)=-nm$, so $\|\mathbf{Y}_n\|=n|m|$. By (2.12) and Theorem 1.5, we have that

$$\left|\frac{\|\mathbf{M}_n\|}{n}-\frac{\|\mathbf{Y}_n\|}{n}\right|\le\frac{\|\widetilde{\mathbf{M}}_n\|}{n}\to 0$$

as $n\to\infty$. This implies that $\|\mathbf{M}_n\|/n\to|m|$ a.s. □

In the context of this paper, the next lemma is very handy for truncation purposes.

Lemma 2.2. *Let $\{X_{ij}: j>i\ge 1\}$ be an infinite triangular array of i.i.d. random variables with $\mathbb{E}X_{12}=0$ and $\operatorname{Var}(X_{12})=\sigma^2$. Let $X_{ji}=X_{ij}$ for $i<j$ and set $X_{ii}=0$ for all $i\ge 1$. Then,*

$$\frac{1}{n^2}\sum_{i=1}^n\left(\sum_{j=1}^n X_{ij}\right)^2\to\sigma^2 \qquad \textit{a.s.}$$

*as $n\to\infty$.*

Proof. Define

$$U_n:=\sum_{i=1}^n\sum_{1\le j<k\le n} X_{ij}X_{ik}. \tag{2.13}$$

Then

$$\frac{1}{n^2}\sum_{i=1}^n\left(\sum_{j=1}^n X_{ij}\right)^2=\frac{1}{n^2}\sum_{i=1}^n\sum_{j=1}^n X_{ij}^2+\frac{2}{n^2}U_n.$$

By the strong law of large numbers, the first term on the right-hand side converges almost surely to $\sigma^2$, so it suffices to show that

$$\frac{U_n}{n^2}\to 0 \qquad \text{a.s.} \tag{2.14}$$

To this end, denote by $\mathcal{F}_k$ the $\sigma$-algebra generated by the random variables $\{X_{ij}, 1 \le i, j \le k\}$. Noting that

$$U_{n+1} - U_n = \sum_{1 \le j < k \le n} X_{(n+1)j} X_{(n+1)k} + \sum_{i=1}^{n} \sum_{j=1}^{n} X_{ij} X_{i(n+1)},$$

it is easy to verify that $\{U_n : n \ge 1\}$ is a martingale for the filtration $\{\mathcal{F}_n : n \ge 1\}$. Further, the $n^2(n-1)/2$ terms in the sum (2.13) are uncorrelated. Indeed, if $i \ne i'$ and $j < k$, $j' < k'$, then $\mathbb{E}(X_{ij} X_{ik} X_{i'j'} X_{i'k'}) = 0$ as at least one of the four variables in this product must be independent of the others. Thus, $\mathbb{E}(U_n^2) \le \sigma^4 n^2 (n-1)/2$ for any $n \ge 2$, and by Doob's submartingale inequality

$$\mathbb{P}\Bigl(\max_{1 \le i \le m^2} |U_i| \ge m^4 \varepsilon\Bigr) \le \frac{\mathbb{E}(U_{m^2}^2)}{m^8 \varepsilon^2} \le \frac{\sigma^4}{m^2 \varepsilon^2}.$$

It follows by the Borel–Cantelli lemma, that almost surely

$$Z_m := m^{-4} \max_{1 \le i \le m^2} |U_i| \to 0,$$

as $m \to \infty$. Since $n^{-2} |U_n| \le (m/(m-1))^4 Z_m$ whenever $(m-1)^2 \le n \le m^2$, $m \ge 2$, we thus get (2.14). □

Let $d_{\mathrm{BL}}$ denote the bounded Lipschitz metric

$$d_{\mathrm{BL}}(\mu, \nu) = \sup\Bigl\{\int f\, d\mu - \int f\, d\nu : \|f\|_\infty + \|f\|_L \le 1\Bigr\}, \tag{2.15}$$

where $\|f\|_\infty = \sup_b |f(x)|$, $\|f\|_L = \sup_{x \ne y} |f(x) - f(y)|/|x-y|$. It is well known (see [8], Section 11.3) that $d_{\mathrm{BL}}$ is a metric for the weak convergence of measures. For the spectral measures of $n \times n$ symmetric real matrices $\mathbf{A}, \mathbf{B}$ we have

$$d_{\mathrm{BL}}(\hat\mu(\mathbf{A}), \hat\mu(\mathbf{B})) \le \sup\Biggl\{\frac{1}{n} \sum_{j=1}^{n} |f(\lambda_j(\mathbf{A})) - f(\lambda_j(\mathbf{B}))| : \|f\|_L \le 1\Biggr\}$$

$$\le \frac{1}{n} \sum_{j=1}^{n} |\lambda_j(\mathbf{A}) - \lambda_j(\mathbf{B})|.$$

By Lidskii's theorem ([13], see also [1], Lemma 2.3)

$$\sum_{j=1}^{n} |\lambda_j(\mathbf{A}) - \lambda_j(\mathbf{B})|^2 \le \operatorname{tr}((\mathbf{B} - \mathbf{A})^2),$$

so

$$d_{\mathrm{BL}}^2(\hat\mu(\mathbf{A}), \hat\mu(\mathbf{B})) \le \frac{1}{n} \operatorname{tr}((\mathbf{B} - \mathbf{A})^2). \tag{2.16}$$

PROOF OF THEOREM 1.4. We use the notation from the proof of Corollary 1.6 and write $\sigma^2 = \mathrm{Var}(X_{11})$. By (2.12) and (2.16) the bounded Lipschitz metric (2.15) satisfies

$$d_{\mathrm{BL}}(\hat{\mu}(\mathbf{M}_n/n), \hat{\mu}(\mathbf{Y}_n/n)) \leq (n^{-3}\,\mathrm{tr}(\widetilde{\mathbf{M}}_n^2))^{1/2}. \tag{2.17}$$

Note that $\{\widetilde{X}_{ij}; 1 \leq i \leq j\}$ are i.i.d. random variables with mean zero and finite variance. By the classical strongumbers and Lemma 2.2

$$n^{-2}\,\mathrm{tr}(\widetilde{\mathbf{M}}_n^2) = \left(\frac{2}{n^2}\sum_{1\leq i<j\leq n}\widetilde{X}_{ij}^2 + \frac{1}{n^2}\sum_{i=1}^{n}\left(\sum_{j\neq i}^{n}\widetilde{X}_{ij}\right)^2\right) \to 2\sigma^2 \qquad \text{a.s.} \tag{2.18}$$

as $n \to \infty$. Recall that all but one of the eigenvalues of $\mathbf{Y}_n$ are $-nm$, hence $\hat{\mu}(\mathbf{Y}_n/n)$ converges weakly to $\delta_{-m}$. Combining this with (2.17) and (2.18), we have that almost surely, $\hat{\mu}(\mathbf{M}_n/n)$ converges weakly to $\delta_{-m}$. $\square$

## 3. The limiting distributions $\gamma_H$, $\gamma_M$ and $\gamma_T$.

3.1. *Moments.* For a probability measure $\gamma$ on $(\mathbb{R}, \mathcal{B})$, denote its moments by

$$m_k(\gamma) = \int x^k \gamma(dx).$$

The probability measures $\gamma_H$, $\gamma_M$ and $\gamma_T$ will be determined from their moments. It turns out that the odd moments are zero, and the even moments are the sums of numbers labeled by the pair partitions of $\{1, \ldots, 2k\}$.

It is convenient to index the pair partitions by the *partition words* $w$; these are words of length $|w| = 2k$ with $k$ pairs of letters such that the first occurrences of each of the $k$ letters are in alphabetic order. In the case $k = 2$ we have $1 \times 3$ such partition words

$$aabb \qquad abba \qquad abab,$$

which correspond to the pair partitions

$$\{1, 2\} \cup \{3, 4\} \qquad \{1, 4\} \cup \{2, 3\} \qquad \{1, 3\} \cup \{2, 4\}$$

of $\{1, 2, 3, 4\}$. Recall that the number of pair partitions of $\{1, \ldots, 2k\}$ is $1 \times 3 \times \cdots \times (2k - 1)$.

DEFINITION 3.1. For a partition word $w$, we define its *height* $h(w)$ as the number of *encapsulated partition subwords*, that is, substrings of the form $xw_1x$, where $x$ is a single letter, and $w_1$ is either a partition word or the empty word.

For example, $h(abcabc) = 0$, $h(\underline{a}bcbc\underline{a}) = h(ab\underline{cc}ab) = 1$, while $h(\underline{aa}\underline{bb}\underline{cc}) = h(\underline{a}\underline{b}\underline{cc}\underline{b}\underline{a}) = 3$ (the encapsulating pairs of letters are underlined).

In the terminology of Bożejko and Speicher [5], $h$ assigns to a pair partition the number of connected blocks which are of cardinality 2. These connected blocks of cardinality 2 are the pairs of letters underlined in the previous examples.

In Proposition A.5 we show that the even moments of the free convolution $\gamma_M$ of the semicircle and standard normal measures are given by

$$m_{2k}(\gamma_M) = \sum_{w\,:\,|w|=2k} 2^{h(w)}. \tag{3.1}$$

For the Toeplitz and Hankel cases, with each partition word $w$ we associate a system of linear equations which determine the cross section of the unit hypercube, and define the corresponding volume $p(w)$. We have to consider these two cases separately.

3.2. *Toeplitz volumes.* Let $w[j]$ denote the letter in position $j$ of the word $w$. For example, if $w = abab$, then $w[1] = a, w[2] = b, w[3] = a, w[4] = b$.

To every partition word $w$ we associate the following system of equations in unknowns $x_0, x_1, \dots, x_{2k}$:

$$\begin{aligned} & x_1 - x_0 + x_{m_1} - x_{m_1-1} = 0, \\ & \qquad \text{if } m_1 > 1 \text{ is such that } w[1] = w[m_1], \\ & x_2 - x_1 + x_{m_2} - x_{m_2-1} = 0, \\ & \qquad \text{if there is } m_2 > 2 \text{ such that } w[2] = w[m_2], \\ & \qquad \vdots \\ & x_i - x_{i-1} + x_{m_i} - x_{m_i-1} = 0, \\ & \qquad \text{if there is } m_i > i \text{ such that } w[i] = w[m_i], \\ & \qquad \vdots \\ & x_{2k-1} - x_{2k-2} + x_{2k} - x_{2k-1} = 0, \\ & \qquad \text{if } w[2k-1] = w[2k]. \end{aligned} \tag{3.2}$$

Although we list $2k-1$ equations, in fact $k-1$ of them are empty. Informally, the left-hand sides of the equations are formed by adding the differences over the same letter when the variables are written in the space "between the letters." For example, writing the variables between the letters of the word $w = ababc..c..$ we get

$${}^{x_0}a^{x_1}b^{x_2}a^{x_3}b^{x_4}c^{x_5}\dots{}^{x_n}c^{x_{n+1}}\dots. \tag{3.3}$$

The corresponding system of equations is

$$\begin{aligned} x_1 - x_0 + x_3 - x_2 &= 0, \\ x_2 - x_1 + x_4 - x_3 &= 0, \\ x_5 - x_4 + x_{n+1} - x_n &= 0, \\ &\vdots. \end{aligned} \tag{3.4}$$

Since in every partition word $w$ of length $2k$ there are exactly $k$ distinct letters, this is the system of $k$ equations in $2k+1$ unknowns. We solve it for the variables that follow the last occurrence of a letter, leaving us with $k+1$ undetermined variables: $x_0$, and the $k$ variables that follow the first occurrence of each letter.

We then require that the dependent variables lie in the interval $I = [0,1]$. This determines a cross section of the cube $I^{k+1}$ in the remaining undetermined $k+1$ coordinates, the volume of which we denote by $p_T(w)$. For example, if $w = abab$, solving the first pair of equations (3.4) for $x_3 = x_0 - x_1 + x_2$, $x_4 = x_0$, defines the solid

$$\{x_0 - x_1 + x_2 \in I\} \cap \{x_0 \in I\} \subset I^3,$$

which has the (Eulerian) volume $p_T(abab) = 4/3! = 2/3$.

We define measure $\gamma_T$ as a symmetric measure with even moments

$$m_{2k}(\gamma_T) = \sum_{w\,:\,|w|=2k} p_T(w). \tag{3.5}$$

From Proposition 4.5 below it follows that (3.5) indeed defines a positive definite sequence of numbers so that these are indeed the even moments of a probability measure. Since $m_{2k}$ is at most the number $(2k-1)!!$ of words of length $2k$, these moments determine the limiting distribution $\gamma_T$ uniquely.

3.3. *Hankel volumes.* We proceed similarly to the Toeplitz case. With each partition word $w$ we associate the following system of equations in unknowns $x_0, x_1, \ldots, x_{2k}$:

$$\begin{aligned} x_1 + x_0 &= x_{m_1} + x_{m_1-1}, \\ &\qquad \text{if } m_1 > 1 \text{ is such that } w[1] = w[m_1], \\ x_2 + x_1 &= x_{m_2} + x_{m_2-1}, \\ &\qquad \text{if there is } m_2 > 2 \text{ such that } w[2] = w[m_2], \\ &\vdots \\ x_i + x_{i-1} &= x_{m_i} + x_{m_i-1}, \end{aligned} \tag{3.6}$$

$$\text{if there is } m_i > i \text{ such that } w[i] = w[m_j],$$

$$\vdots$$

$$x_{2k-1} + x_{2k-2} = x_{2k} + x_{2k-1},$$
$$\text{if } w[2k-1] = w[2k].$$

Informally, the equations are formed by equating the sums of the variables at the same letter. For example, the word *abab* with the variables written as in (3.3) gives rise to the system of equations

$$\begin{aligned} x_1 + x_0 &= x_3 + x_2, \\ x_2 + x_1 &= x_4 + x_3. \end{aligned} \tag{3.7}$$

As in the Toeplitz case, since there are exactly $k$ distinct letters in the word, this is the system of $k$ equations in $2k+1$ unknowns. We solve it for the variables that precede the first occurrence of a letter, leaving us with $k$ undetermined variables $\dots, x_{\alpha_1}, \dots, x_{\alpha_k} = x_{2k-1}$ that precede the second occurrence of each letter, and with the $(k+1)$st undetermined variable $x_{2k}$. We add to the system (3.6) one more equation:

$$x_0 = x_{2k}.$$

As previously, we require that the dependent variables are in the interval $I = [0, 1]$. This determines a cross section of the cube $I^{k+1}$ in the remaining $k+1$ coordinates with the volume which we denote by $p_H(w)$.

Due to the additional constraint $x_{2k} = x_0$, this volume might be zero. For example, (3.7) has solutions $x_0 = 2x_2 - x_4, x_1 = x_3 - x_2 + x_4$ with undetermined variables $x_2, x_3, x_4$. Equation $x_0 = x_4$ gives additional relation $x_4 = x_2$, and reduces the dimension of the solid $\{2x_2 - x_4 \in I\} \cap \{x_3 - x_2 + x_4 \in I\} \cap \{x_4 = x_2\} \subset I^3$ to 2. Thus the corresponding volume is $p_H(abab) = 0$.

We define measure $\gamma_H$ as a symmetric measure with even moments

$$m_{2k}(\gamma_H) = \sum_{w\,:\,|w|=2k} p_H(w). \tag{3.8}$$

From Proposition 4.7 below it follows that (3.8) indeed defines a positive definite sequence of numbers so that these are indeed the even moments of a probability measure. Since $m_{2k}$ is at most the number $(2k-1)!!$ of words of length $2k$, these moments determine the limiting distribution $\gamma_H$ uniquely.

3.4. *Relation to Eulerian numbers.* The Eulerian numbers $A_{n,m}$ are often defined by their generating function or by the combinatorial description as the number of permutations $\sigma$ of $\{1, \dots, n\}$ with $\sigma_i > \sigma_{i-1}$ for exactly $m$ choices of $i = 1, 2, \dots, n$ (taking $\sigma_0 = 0$). The geometric interpretation is that $A_{n,m}/n!$ is the volume of a solid cut out of the cube $I^n$ by the set $\{x_1+\cdots+x_n \in [m-1,m]\}$; see [23]. Converting any $m-1$ of the coordinates $x$ to $1-x$, we get that $A_{n,m}/n!$ is the volume of a solid cut out of the cube $I^n$ by the set

$$\{(x_1,\dots,x_n)\in\mathbb{R}^n : x_1+x_2+\cdots+x_{n-m}-(x_{n-m+1}+\cdots+x_n)\in I\}.$$

The solids we encountered in the formula for the $2k$th moments are the intersections of solids of this latter form, with odd values of $n$, each having $m=(n-1)/2$, and with various subsets of the coordinates entering the expression.

Another interesting representation is

$$\begin{aligned}&\mathrm{Vol}(\{(x_1,\dots,x_n)\in I^n : x_1+x_2+\cdots+x_{n-m}-(x_{n-m+1}+\cdots+x_n)\in I\})\\&\quad=\frac{2}{\pi}\int_0^\infty\left(\frac{\sin t}{t}\right)^{n+1}\cos((n+1-2m)t)\,dt.\end{aligned}$$

This follows from the integral representation of Eulerian numbers in Nicolas [16].

REMARK 3.1. One can verify that the probabilities $p_T(w)$ and $p_H(w)$ are rational numbers, and hence so are $m_{2k}(\gamma_T)$ and $m_{2k}(\gamma_H)$, defined by formulas (3.5) and (3.8) (for details, cf. [6]).

## 4. Proofs of Theorems 1.1, 1.2 and 1.3.

4.1. *Truncation and centering.* We first reduce Theorems 1.1, 1.2 and 1.3 to the case of bounded i.i.d. random variables, and in case of Theorems 1.1 and 1.2, also allow for centering of these variables.

PROPOSITION 4.1. (i) *If Theorem* 1.1 *or Theorem* 1.2 *holds true for all bounded independent i.i.d. sequences* $\{X_j\}$ *with mean zero and variance* 1, *then it holds true for all square-integrable i.i.d. sequences* $\{X_j\}$ *with variance* 1.

(ii) *If Theorem* 1.3 *holds true for all bounded independent i.i.d. collections* $\{X_{ij}\}$ *with mean zero and variance* 1, *then it holds true for all square-integrable i.i.d. collections* $\{X_{ij}\}$ *with mean zero and variance* 1.

PROOF. Without loss of generality, we may assume that $\mathbb{E}(X_1)=0$ in Theorems 1.1 and 1.2. Indeed, from the rank inequality ([1], Lemma 2.2) it follows that subtracting a rank-1 matrix of the means $\mathbb{E}(X_1)$ from matrices $\mathbf{T}_n$ and $\mathbf{H}_n$ does not affect the asymptotic distribution of the eigenvalues.

For a fixed $u>0$, denote

$$m(u)=\mathbb{E}X_1 I_{\{|X_1|>u\}}$$

and let

$$\sigma^2(u) = \mathbb{E}X_1^2 I_{\{|X_1|\le u\}} - m^2(u).$$

Clearly, $\sigma^2(u) \le 1$ and since $\mathbb{E}(X_1) = 0$, $\mathbb{E}(X_1^2) = 1$, we have $m(u) \to 0$ and $\sigma(u) \to 1$ as $u \to \infty$.

Let

$$\widetilde{X}_1 = X_1 I_{\{|X_1|>u\}} - m(u).$$

Notice that $\sigma^2(u) = \mathbb{E}(X_1 - \widetilde{X}_1)^2$, therefore the bounded random variable

$$X_1' = \frac{X_1 - \widetilde{X}_1}{\sigma(u)}$$

has mean zero and variance 1. Denote by $\mathbf{T}_n', \mathbf{H}_n'$ the corresponding Toeplitz and Hankel matrices constructed from the independent bounded random variables

$$X_j' := \frac{X_j - \widetilde{X}_j}{\sigma(u)}$$

distributed as $X_1'$. By the triangle inequality for $d_{\mathrm{BL}}(\cdot,\cdot)$ and (2.16),

$$\begin{aligned} &d_{\mathrm{BL}}^2(\hat{\mu}(\mathbf{T}_n/\sqrt{n}), \hat{\mu}(\mathbf{T}_n'/\sqrt{n})) \\ &\quad \le 2d_{\mathrm{BL}}^2(\hat{\mu}(\mathbf{T}_n/\sqrt{n}), \hat{\mu}(\sigma(u)\mathbf{T}_n'/\sqrt{n})) + 2d_{\mathrm{BL}}^2(\hat{\mu}(\mathbf{T}_n'/\sqrt{n}), \hat{\mu}(\sigma(u)\mathbf{T}_n'/\sqrt{n})) \\ &\quad \le \frac{2}{n^2}\operatorname{tr}((\mathbf{T}_n - \sigma(u)\mathbf{T}_n')^2) + \frac{2}{n^2}(1-\sigma(u))^2 \operatorname{tr}((\mathbf{T}_n')^2). \end{aligned}$$

It is easy to verify that $\mathbb{E}(\widetilde{X}_1^2) = 1 - \sigma^2(u) - 2m(u)^2$ and that with probability 1

$$\frac{1}{n^2}\operatorname{tr}((\mathbf{T}_n - \sigma(u)\mathbf{T}_n')^2) = \frac{1}{n}\widetilde{X}_0^2 + \frac{2}{n}\sum_{j=1}^{n}\left(1 - \frac{j}{n}\right)\widetilde{X}_j^2 \to \mathbb{E}(\widetilde{X}_1^2), \tag{4.1}$$

as $n \to \infty$ (e.g., sandwiching the coefficients $j/n$ between the piecewise constant $\ell^{-1}\lfloor \ell j/n \rfloor$ and $\ell^{-1}\lceil \ell j/n \rceil$ allows for applying the strong law of large numbers, with the resulting nonrandom bounds converging to $\mathbb{E}(\widetilde{X}_1^2)$ as $\ell \to \infty$). Similarly,

$$\frac{1}{n^2}\operatorname{tr}((\mathbf{T}_n')^2) = \frac{1}{n}(X_0')^2 + \frac{2}{n}\sum_{j=1}^{n}\left(1 - \frac{j}{n}\right)(X_j')^2 \to \mathbb{E}((X_1')^2). \tag{4.2}$$

For large $u$, both $m(u)$ and $1-\sigma(u)$ are arbitrarily small. So, in view of (4.1) and (4.2), with probability 1 the limiting distance in the bounded Lipschitz metric $d_{\mathrm{BL}}$ between $\hat{\mu}(\mathbf{T}_n/\sqrt{n})$ and $\hat{\mu}(\mathbf{T}_n'/\sqrt{n})$ is arbitrarily small, for all $u$ sufficiently large. Thus, if the conclusion of Theorem 1.1 holds true for all sequences of independent bounded random variables $\{X'_j\}$, with the same limiting distribution $\gamma_T$, then $\hat{\mu}(\mathbf{T}_n/\sqrt{n})$ must have the same weak limit with probability 1.

Similarly, we have

$$d^2_{\mathrm{BL}}(\hat{\mu}(\mathbf{H}_n/\sqrt{n}), \hat{\mu}(\mathbf{H}'_n/\sqrt{n})) \\ \leq \frac{2}{n^2}\operatorname{tr}((\mathbf{H}_n - \sigma(u)\mathbf{H}'_n)^2) + \frac{2}{n^2}(1-\sigma(u))^2 \operatorname{tr}((\mathbf{H}'_n)^2).$$

By the same argument as before, with probability 1

$$\frac{1}{n^2}\operatorname{tr}((\mathbf{H}_n - \sigma(u)\mathbf{H}'_n)^2) = \frac{1}{n}\sum_{j=0}^{2n}\left(1 - \frac{|j-n|}{n}\right)\widetilde{X}_j^2 \to \mathbb{E}(\widetilde{X}_1^2),$$

and $n^{-2}\operatorname{tr}((\mathbf{H}'_n)^2) \to \mathbb{E}((X'_1)^2)$. Therefore, with probability 1 the limiting $d_{\mathrm{BL}}$-distance between $\hat{\mu}(\mathbf{H}_n/\sqrt{n})$ and $\hat{\mu}(\mathbf{H}'_n/\sqrt{n})$ is arbitrarily small for large enough $u$.

Similarly, denoting by $\widetilde{\mathbf{M}}_n, \mathbf{M}'_n$ the corresponding Markov matrices constructed from the independent bounded random variables $\widetilde{X}_{ij}$ and $X'_{ij} := \frac{X_{ij} - \widetilde{X}_{ij}}{\sigma(u)}$, we have

$$d^2_{\mathrm{BL}}(\hat{\mu}(\mathbf{M}_n/\sqrt{n}), \hat{\mu}(\mathbf{M}'_n/\sqrt{n})) \leq \frac{2}{n^2}\operatorname{tr}(\widetilde{\mathbf{M}}_n^2) + \frac{2}{n^2}(1-\sigma(u))^2 \operatorname{tr}((\mathbf{M}'_n)^2).$$

By (2.18), with probability 1, $n^{-2}\operatorname{tr}((\mathbf{M}'_n)^2) \to 2$ and $n^{-2}\operatorname{tr}(\widetilde{\mathbf{M}}_n^2) \to 2\mathbb{E}(\widetilde{X}_{12}^2)$. Therefore, with probability 1, the limiting $d_{\mathrm{BL}}$-distance between $\hat{\mu}(\mathbf{M}_n/\sqrt{n})$ and $\hat{\mu}(\mathbf{M}'_n/\sqrt{n})$ is arbitrarily small for large enough $u$. $\square$

4.2. *Combinatorics for Hankel and Toeplitz cases.* For $k, n \in \mathbb{N}$, consider circuits in $\{1, \dots, n\}$ of length $L(\pi) = k$, that is, mappings $\pi : \{0, 1, \dots, k\} \to \{1, 2, \dots, n\}$, such that $\pi(0) = \pi(k)$.

Let $s : \mathbb{N}^2 \to \mathbb{N}$ be one of the following two functions: $s_T(x, y) = |x - y|$, or $s_H(x, y) = x + y$. We will use $s$ to match (i.e., pair) the edges $(\pi(i-1), \pi(i))$ of a circuit $\pi$. The main property of the symmetric function $s$ is that for a fixed value of $s(m, n)$, every initial point $m$ of an edge determines uniquely a finite number (here, at most 2) of the other end-points: if $k, m \in \mathbb{N}$, then

$$\#\{y \in \mathbb{N} : s(m, y) = k\} \leq 2. \tag{4.3}$$

For a fixed $s$ as above, we will say that circuit $\pi$ is $s$-matched, or has self-matched edges, if for every $1 \leq i \leq L(\pi)$ there is $j \neq i$ such that $s(\pi(i-1), \pi(i)) = s(\pi(j-1), \pi(j))$.

We will say that a circuit $\pi$ has an edge of order 3, if there are at least three different edges in $\pi$ with the same $s$-value.

The following proposition says that generically self-matched circuits have only pair-matches.

Proposition 4.2. *Fix $r \in \mathbb{N}$. Let $N$ denote the number of $s$-matched circuits in $\{1,\dots,n\}$ of length $r$ with at least one edge of order 3. Then there is a constant $C_r$ such that*

$$N \le C_r n^{\lfloor (r+1)/2 \rfloor}.$$

*In particular, as $n \to \infty$ we have $\frac{N}{n^{1+r/2}} \to 0$.*

Proof. Either $r = 2k$ is an even number, or $r = 2k-1$ is an odd number. In both cases, if an $s$-matched circuit has an edge of order 3, then the total number of distinct $s$-values

$$\{s(\pi(i-1), \pi(i)) : 1 \le i \le L(\pi)\}$$

is at most $k-1$. We can think of constructing each such circuit from the left to the right. First, we choose the locations for the $s$-matches along $\{1,\dots,r\}$. This can be done in at most $r!$ ways. Once these locations are fixed, we proceed along the circuit. There are $n$ possible choices for the initial point $\pi(0)$. There are at most $n$ choices for each new $s$-value, and there are at most two ways to complete the edge for each repeat of the already encountered $s$-value. Therefore there are at most $r! \times n \times n^{k-1} 2^{r+1-k} \le C_r n^k$ such circuits. $\square$

We say that a set of circuits $\pi_1, \pi_2, \pi_3, \pi_4$ is matched if each edge of any one of these circuits is either self-matched, that is, there is another edge of the same circuit with equal $s$-value, or is cross-matched, that is, there is an edge of the other circuit with the same $s$-value (or both).

The following bound will be used to prove almost sure convergence of moments.

Proposition 4.3. *Fix $r \in \mathbb{N}$. Let $N$ denote the number of matched quadruples of circuits in $\{1,\dots,n\}$ of length $r$ such that none of them is self-matched. Then there is a constant $C_r$ such that*

$$N \le C_r n^{2r+2}.$$

Proof. First observe that there are at most $2r$ distinct $s$-values in the $4r$ edges of matched quadruples of circuits of length $r$. Further, the number of quadruples of such circuits for which there are exactly $u$ distinct $s$-values is at most $C_{r,u} n^{u+4}$. Indeed, order the edges $(\pi_j(i-1), \pi_j(i))$, of such quadruples starting at $j = 1$, $i = 1$, then $i = 2,\dots,r$, followed by $j = 2$, $i = 1$ and then $i = 2,\dots,r$, and so on. There are at most $u^{4r}$ possible allocations of the distinct $s$-values to these $4r$ edges, at most $n^4$ choices for the starting points $\pi_1(0)$, $\pi_2(0)$, $\pi_3(0)$ and $\pi_4(0)$ of the circuits and at most $n^u$ for the values of $\pi_j(i)$ at those $(j,i)$ for which $(\pi_j(i-1), \pi_j(i))$ is the leftmost occurrence of one of the distinct $s$-values. Once these choices are made, we proceed to sequentially determine the mapping $\pi_1(i)$ from $i=0$ to $i=r$, followed by the mappings $\pi_2, \pi_3, \pi_4$, noting that by (4.3) at most $2^{4r-u-4}$ quadruples can be produced per such choice.

Recall that the number of possible partitions $\mathcal{P}$ of the $4r$ edges of our quadruple of circuits into $|\mathcal{P}|$ distinct groups of $s$-matching edges, with at least two edges in each group, is independent of $n$. Thus, by the preceding bound it suffices to show that for each partition $\mathcal{P}$ with $|\mathcal{P}| \in \{2r-1, 2r\}$ such that *each circuit shares at least one $s$-value with some other circuit*, there correspond at most $Cn^{2r+2}$ matched quadruples of circuits in $\{1, \dots, n\}$. To this end, note that $|\mathcal{P}| = 2r$ implies that each $s$-value is shared by exactly two edges, while when $|\mathcal{P}| = 2r-1$ we also have either two $s$-values shared by three edges each or one $s$-value shared by four edges (but not both).

Fixing hereafter a specific partition $\mathcal{P}$ of this type, it is not hard to check that upon re-ordering our four circuits we have an $s$-value that is assigned to exactly one edge of the circuit $\pi_1$, denoted hereafter $(\pi_1(i_*-1), \pi_1(i_*))$, and in case $|\mathcal{P}| = 2r$, we also have another $s$-value that does not appear in $\pi_1$ and is assigned to exactly one edge of $\pi_2$, denoted hereafter $(\pi_2(j_*-1), \pi_2(j_*))$. (Though this property may not hold for all ordering of the four circuits, an inspection of all possible graphs of cross-matches shows that it must hold for some order.)

We are now ready to improve our counting bound for the case of $|\mathcal{P}| = 2r-1$, by the following dynamic construction of $\pi_1$:

First choose one of the $n$ possible values for the initial value $\pi_1(0)$, and continue filling in the values of $\pi_1(i)$, $i = 1, 2, \dots, i_*-1$. Then, starting at $\pi_1(r) = \pi_1(0)$, sequentially choose the values of $\pi_1(r-1), \pi_1(r-2), \dots, \pi_1(i^*)$, thus completing the entire circuit $\pi_1$. This is done in accordance with the $s$-matches determined by $\mathcal{P}$, so there are $n$ ways to complete an edge that has no $s$-match among the edges already constructed, while by (4.3) if an edge is matching one of the edges already available, then it can be completed in at most two ways. Since this procedure determines uniquely the edge $(\pi_1(i_*-1), \pi_1(i_*))$ and hence the $s$-value assigned to it, it reduces to $2r-2$ the number of $s$-matches that can each independently assume $O(n)$ values. Consequently, the number of quadruples of circuits corresponding to $\mathcal{P}$ is at most $Cn^{2r+2}$.

In case $|\mathcal{P}| = 2r$, we first construct $\pi_1$ by the preceding dynamic construction while determining the $s$-value for the edge $(\pi_1(i_*-1), \pi_1(i_*))$ out of the circuit condition for $\pi_1$. Then, we repeat the dynamic construction for $\pi_2$, keeping it in accordance with the $s$-values determined already by edges of $\pi_1$ and uniquely determining the edge $(\pi_2(j_*-1), \pi_2(j_*))$ and hence the $s$-value assigned to it, by the circuit condition for $\pi_2$. Thus, we have again reduced the total number of $s$-matches that can each independently assume $O(n)$ values to $2r-2$, and consequently, the number of quadruples of circuits corresponding to $\mathcal{P}$ is again at most $Cn^{2r+2}$. $\square$

The next result deals only with the slope matching function $s_T(x,y) = |x-y|$.

PROPOSITION 4.4. *Fix $k \in \mathbb{N}$. Let $N$ be the number of $s_T$-matched circuits $\pi$ in $\{1,\dots,n\}$ of length $2k$ with at least one pair of $s_T$-matched edges $(\pi(i-1),\pi(i))$ and $(\pi(j-1),\pi(j))$ such that $\pi(i)-\pi(i-1)+\pi(j)-\pi(j-1) \neq 0$. Then, as $n \to \infty$ we have*

$$n^{-(k+1)}N \to 0.$$

PROOF. By Proposition 4.2, we may and shall consider throughout path $\pi$ in $\{1,\dots,n\}$ of length $2k$ for which the absolute values of the slopes $\pi(i)-\pi(i-1)$ take exactly $k$ distinct nonzero values and, for $\pi$ to be a circuit, the sum of all $2k$ slopes is zero. Let $\mathcal{P}$ denote a partition of the $2k$ slopes to $s_T$-matching pairs, indicating also whether each slope is negative or positive, with $m(\mathcal{P})$ denoting the number of such pairs for which both slopes are positive. Observe that if under $\mathcal{P}$ both slopes of some $s_T$-matching pair are negative, then necessarily $m(\mathcal{P}) \geq 1$, for otherwise the sum of all slopes will not be zero for any path corresponding to $\mathcal{P}$. Thus, it suffices to show that at most $n^k$ circuits $\pi$ correspond to each $\mathcal{P}$ with $m = m(\mathcal{P}) \geq 1$. Indeed, fixing such $\mathcal{P}$, there are at most $n$ ways to choose $\pi(0)$ and $n^{k-m}$ ways to choose the $k-m$ pairs of slopes for which at least one slope in each pair is negative. The remaining $m$ pairs of $s_T$-matching positive slopes are to be chosen among $\{1,\dots,n\}$ subject to a specified sum (due to the circuit condition). Since there are at most $n^{m-1}$ ways for doing so, the proof is complete. $\square$

4.3. *Moments of the average spectral measure.*

PROPOSITION 4.5. *Suppose $\{X_j\}$ is a sequence of bounded i.i.d. random variables such that $\mathbb{E}(X_1)=0, \mathbb{E}(X_1^2)=1$. Then for $k \in \mathbb{N}$*

$$\lim_{n\to\infty} \frac{1}{n^{k+1}} \mathbb{E}\,\mathrm{tr}(\mathbf{T}_n^{2k}) = \sum_{w\,:\,|w|=2k} p_T(w) \tag{4.4}$$

*and*

$$\lim_{n\to\infty} \frac{1}{n^{k+1/2}} \mathbb{E}\,\mathrm{tr}(\mathbf{T}_n^{2k-1}) = 0. \tag{4.5}$$

PROOF. For a circuit $\pi:\{0,1,\dots,r\}\to\{1,2,\dots,n\}$ write

$$\mathbf{X}_\pi=\prod_{i=1}^{r}X_{\pi(i)-\pi(i-1)}. \tag{4.6}$$

Then

$$\mathbb{E}\,\mathrm{tr}(\mathbf{T}_n^r)=\sum_{\pi}\mathbb{E}\mathbf{X}_\pi, \tag{4.7}$$

where the sum is over all circuits in $\{1,\dots,n\}$ of length $r$.

By Hölder's inequality, for any finite set $\Pi$ of circuits of length $r$

$$\left|\sum_{\pi\in\Pi}\mathbb{E}\mathbf{X}_\pi\right|\le\mathbb{E}(|X|^r)\#\Pi. \tag{4.8}$$

Since $|X|^r$ is bounded, we can use the bound (4.8) to discard the "nongeneric" circuits from the sum in (4.7). To this end, note that since the random variables $\{X_j\}$ are independent and have mean zero, the term $\mathbb{E}\mathbf{X}_\pi$ vanishes for every circuit $\pi$ with at least one unpaired $X_j$. Since $\mathbf{T}_n$ is a symmetric matrix, by (4.6) paired variables correspond to the slopes of the circuit $\pi$ which are equal in absolute value. Hence, the only circuits that make a nonzero contribution to (4.7) are those with matched absolute values of the slopes. This fits the formalism of Section 4.2 with the matching function $s_T(x,y)=|x-y|$.

If $r=2k-1>0$ is odd, then each $s_T$-matched circuit $\pi$ of length $r$ must have an edge of order 3. From (4.8) and Proposition 4.2 we get $|\mathbb{E}\,\mathrm{tr}(\mathbf{T}_n^{2k-1})|\le Cn^k$, proving (4.5).

When $r=2k$ is an even number, let $\Pi$ be the set of all circuits $\pi:\{0,1,\dots,2k\}\to\{1,\dots,n\}$ with the set of slopes $\{\pi(i)-\pi(i-1):i=1,\dots,2k\}$ consisting of $k$ distinct nonnegative integers $s_1,\dots,s_k$ and their counterparts $-s_1,\dots,-s_k$. From (4.8) and Proposition 4.4 it follows that

$$\lim_{n\to\infty}\frac{1}{n^{k+1}}\left|\mathbb{E}\,\mathrm{tr}(\mathbf{T}_n^r)-\sum_{\pi\in\Pi}\mathbb{E}\mathbf{X}_\pi\right|=0.$$

Moreover, for every circuit $\pi\in\Pi$, if $X_j$ enters the product $\mathbf{X}_\pi$, then it occurs in it exactly twice, resulting with $\mathbb{E}\mathbf{X}_\pi=1$, and consequently with $\sum_{\pi\in\Pi}\mathbb{E}\mathbf{X}_\pi=\#\Pi$. Therefore, the following lemma completes the proof of (4.4), and with it, that of Proposition 4.5. $\square$

LEMMA 4.6.

$$\lim_{n\to\infty}\frac{1}{n^{k+1}}\#\Pi=\sum_{w}p_T(w),$$

*where the sum is over the finite set of partition words $w$ of length $2k$.*

PROOF. The circuits in $\Pi$ can be labeled by the partition words $w$ of length $2k$ which list the positions of the pairs of $s_T$-matches along $\{1,\dots,2k\}$. This generates the partition $\Pi = \bigcup_w \Pi(w)$ into the corresponding equivalence classes.

To every such partition word $w$ we can assign $n^{k+1}$ paths $\pi(i) = x_i$, $i = 0,\dots,2k$, obtained by solving the system of equations (3.2), with values $1,2,\dots,n$ for each of the $k+1$ undetermined variables, and the remaining $k$ values computed from the equations [which represent the relevant $s_T$-matches for any $\pi \in \Pi(w)$]. Some of these paths will fail to be in the admissible range $\{1,\dots,n\}$. Let $p_n(w)$ be the fraction of the $n^{k+1}$ paths that stay within the admissible range $\{1,\dots,n\}$, noting that by Proposition 4.2, $p_n(w) - n^{-(k+1)}\#\Pi(w) \to 0$.

Interpreting the undetermined variables $x_j$ as the discrete uniform independent random variables with values $\{1,2,\dots,n\}$, $p_n(w)$ becomes the probability that the computed values stay within the prescribed range. As $n \to \infty$, the $k+1$ undetermined variables $x_j/n$ converge in law to independent uniform $U[0,1]$ random variables $U_j$. Since $p_n(w)$ is the probability of the (independent of $n$) event $A_w$ that the solution of (3.2) starting with $x_j/n \in \{1/n, 2/n,\dots,1\}$ has all the dependent variables in $(0,1]$, it follows that $p_n(w)$ converges to $p_T(w)$, the probability of the event $A_w$ that the corresponding sums of independent uniform $U[0,1]$ random variables take their values in the interval $[0,1]$. $\square$

Next we give the Hankel version of Proposition 4.5.

PROPOSITION 4.7. *Let $\{X_j\}$ be a sequence of bounded i.i.d. random variables such that $\mathbb{E}(X_1) = 0, \mathbb{E}(X_1^2) = 1$. For $k \in \mathbb{N}$,*

$$\lim_{n\to\infty} \frac{1}{n^{k+1}} \mathbb{E}\,\mathrm{tr}(\mathbf{H}_n^{2k}) = \sum_{w\,:\,|w|=2k} p_H(w) \tag{4.9}$$

*and*

$$\lim_{n\to\infty} \frac{1}{n^{k+1/2}} \mathbb{E}\,\mathrm{tr}(\mathbf{H}_n^{2k-1}) = 0. \tag{4.10}$$

PROOF. We mimic the procedure for the Toeplitz case. For a circuit $\pi : \{0,1,\dots,r\} \to \{1,2,\dots,n\}$ write

$$\mathbf{X}_\pi = \prod_{i=1}^{r} X_{\pi(i)+\pi(i-1)}. \tag{4.11}$$

As previously,

$$\mathbb{E}\,\mathrm{tr}(\mathbf{H}_n^r) = \sum_\pi \mathbb{E}\mathbf{X}_\pi, \tag{4.12}$$

where the sum is over all circuits in $\{1,\dots,n\}$ of length $r$, and by Hölder's inequality, we again have the bound (4.8), which for bounded $|X|^r$ we use to discard the "nongeneric" circuits from the sum in (4.12). To this end, with the random variables $X_j$ independent and of mean zero, the term $\mathbb{E}\mathbf{X}_\pi$ vanishes for every circuit $\pi$ with at least one unpaired $X_j$. By (4.11), in the current setting paired variables correspond to an $s_H$-matching in the circuit $\pi$. Hence, only $s_H$-matched circuits (in the formalism of Section 4.2) can make a nonzero contribution to (4.12).

If $r=2k-1>0$ is odd, then each $s_H$-matched circuit $\pi$ of length $r$ must have an edge of order 3. From (4.8) and Proposition 4.2 we get $|\mathbb{E}\operatorname{tr}(\mathbf{H}_n^{2k-1})|\le Cn^k$, proving (4.10).

When $r=2k$ is an even number, let $\Pi$ be the set of all circuits $\pi:\{0,1,\dots,2k\}\to\{1,\dots,n\}$ with the $s_H$-values consisting of $k$ distinct numbers. Recall that $\mathbb{E}\mathbf{X}_\pi=1$ for any $\pi\in\Pi$ [see (4.11)]. Further, with any $s_H$-matched circuit not in $\Pi$ having an edge of order 3, it follows from (4.8) and Proposition 4.2 that

$$\lim_{n\to\infty}\frac{1}{n^{k+1}}|\mathbb{E}\operatorname{tr}(\mathbf{H}_n^r)-\#\Pi|=0.$$

Therefore, the following lemma completes the proof of (4.9), and with it, that of Proposition 4.7. $\square$

Lemma 4.8.

$$\lim_{n\to\infty}\frac{1}{n^{k+1}}\#\Pi=\sum_{w:|w|=2k}p_H(w).$$

Proof. Similarly to the proof of Lemma 4.6, label the circuits in $\Pi$ by the partition words $w$ which list the positions of the pairs of $s_H$-matches along $\{1,\dots,2k\}$, with the corresponding partition $\Pi=\bigcup_w\Pi(w)$ into equivalence classes. To every such partition word $w$ we can assign $n^{k+1}$ paths $\pi(i)=x_i$, $i=0,\dots,2k$, obtained by solving the system of equations (3.6), with values $1,2,\dots,n$ for each of the $k+1$ undetermined variables, and the remaining $k$ values computed from the equations. Some of these paths will fail to be a circuit, and some will fail to stay in the admissible range $\{1,\dots,n\}$. Let $p_n(w)$ denote the fraction of the paths that stay within the admissible range $\{1,\dots,n\}$ and are circuits, noting that $p_n(w)-n^{-(k+1)}\#\Pi(w)\to 0$ by Proposition 4.2. Thus, $p_n(w)$ is the probability of the event $A_w$ that the solution of (3.6) starting with the undetermined variables $x_j$ that are independent discrete uniform random variables on the set $\{1/n,2/n,\dots,1\}$, stays within $(0,1]$ and satisfies the additional condition $x_0=x_{2k}$. It follows that as $n\to\infty$, the probabilities $p_n(w)$ converge to $p_H(w)$, the probability of the event $A_w$ with the undetermined variables now being independent and uniformly distributed on $[0,1]$. $\square$

4.4. *Concentration of moments of the spectral measure.*

PROPOSITION 4.9. *Let $\{X_j\}$ be a sequence of bounded i.i.d. random variables such that $\mathbb{E}(X_1)=0$ and $\mathbb{E}(X_1^2)=1$. Fix $r\in\mathbb{N}$. Then there is $C_r<\infty$ such that for all $n\in\mathbb{N}$ we have*

$$\mathbb{E}[(\operatorname{tr}(\mathbf{T}_n^r)-\mathbb{E}\operatorname{tr}(\mathbf{T}_n^r))^4]\le C_r n^{2r+2} \quad \text{and} \quad \mathbb{E}[(\operatorname{tr}(\mathbf{H}_n^r)-\mathbb{E}\operatorname{tr}(\mathbf{H}_n^r))^4]\le C_r n^{2r+2}.$$

PROOF. The argument again relies on the enumeration of paths. Since both proofs are very similar, we analyze only the Hankel case.

Using the circuit notation of (4.11) we have that

$$\mathbb{E}[(\operatorname{tr}(\mathbf{H}_n^r)-\mathbb{E}\operatorname{tr}(\mathbf{H}_n^r))^4]=\sum_{\pi_1,\pi_2,\pi_3,\pi_4}\mathbb{E}\left[\prod_{j=1}^4(\mathbf{X}_{\pi_j}-\mathbb{E}(\mathbf{X}_{\pi_j}))\right], \tag{4.13}$$

where the sum is taken over all circuits $\pi_j$, $j=1,\dots,4$ on $\{1,\dots,n\}$ of length $r$ each. With the random variables $X_j$ independent and of mean zero, any circuit $\pi_k$ which is not matched together with the remaining three circuits has $\mathbb{E}(\mathbf{X}_{\pi_k})=0$ and

$$\mathbb{E}\left[\prod_{j=1}^4(\mathbf{X}_{\pi_j}-\mathbb{E}(\mathbf{X}_{\pi_j}))\right]=\mathbb{E}\left[\mathbf{X}_{\pi_k}\prod_{j\neq k}(\mathbf{X}_{\pi_j}-\mathbb{E}(\mathbf{X}_{\pi_j}))\right]=0.$$

Further, if one of the circuits, say $\pi_1$, is only self-matched, that is, has no cross-matched edge, then obviously

$$\mathbb{E}\left[\prod_{j=1}^4(\mathbf{X}_{\pi_j}-\mathbb{E}(\mathbf{X}_{\pi_j}))\right]=\mathbb{E}[\mathbf{X}_{\pi_1}-\mathbb{E}(\mathbf{X}_{\pi_j})]\mathbb{E}\left[\prod_{j=2}^4(\mathbf{X}_{\pi_j}-\mathbb{E}(\mathbf{X}_{\pi_j}))\right]=0.$$

Therefore, it suffices to take the sum in (4.13) over all $s_H$-matched quadruples of circuits on $\{1,\dots,n\}$, such that none of them is self-matched. By Proposition 4.3, there are at most $C_r n^{2r+2}$ such quadruples of circuits, and with $|X|$ (hence $|\mathbf{X}_\pi|$) bounded, this completes the proof. □

4.5. *Proofs of the Hankel and Toeplitz cases.*

PROOF OF THEOREM 1.1. Proposition 4.1(i) implies that without loss of generality we may assume that the random variables $\{X_j\}$ are centered and bounded.

By Proposition 4.5 the odd moments of the average measure $\mathbb{E}(\hat{\mu}(\mathbf{T}_n/\sqrt{n}\,))$ converge to 0, and the even moments converge to $m_{2k}$ of (3.5). By Chebyshev's inequality we have from Proposition 4.9 that for any $\delta>0$ and $k,n\in\mathbb{N}$,

$$\mathbb{P}\left[\left|\int x^k\,d\hat{\mu}(\mathbf{T}_n/\sqrt{n}\,)-\int x^k\,d\mathbb{E}(\hat{\mu}(\mathbf{T}_n/\sqrt{n}\,))\right|>\delta\right]\le C_k\delta^{-4}n^{-2}.$$

Thus, by the Borel–Cantelli lemma, with probability 1 $\int x^k\, d\hat{\mu}(\mathbf{T}_n/\sqrt{n}\,) \to \int x^k\, d\gamma_T$ as $n\to\infty$, for every $k\in\mathbb{N}$. In particular, with probability 1, the random measures $\{\hat{\mu}(\mathbf{T}_n/\sqrt{n}\,)\}$ are tight, and since the moments determine $\gamma_T$ uniquely, we have the weak convergence of $\hat{\mu}(\mathbf{T}_n/\sqrt{n}\,)$ to $\gamma_T$.

Since the moments do not depend on the distribution of the i.i.d. sequence $\{X_j\}$, the limiting distribution $\gamma_T$ does not depend on the distribution of $X$ either, and is symmetric as all its odd moments are zero. By Proposition A.1, it has unbounded support. $\square$

Proof of Theorem 1.2. We follow the same line of reasoning as in the proof of Theorem 1.1, starting by assuming without loss of generality that $\{X_j\}$ is a sequence of centered and bounded random variables, in view of Proposition 4.1(i). Then, by Proposition 4.7, as $n\to\infty$ the odd moments of the average measure $\mathbb{E}(\hat{\mu}(\mathbf{H}_n/\sqrt{n}\,))$ converge to 0, and the even moments converge to $m_{2k}$ of (3.8), whereas from Proposition 4.9 we conclude that with probability 1 the same applies to the moments of $\hat{\mu}(\mathbf{H}_n/\sqrt{n}\,)$. The almost sure convergence $\int x^k\, d\hat{\mu}(\mathbf{H}_n/\sqrt{n}\,) \to \int x^k\, d\gamma_H$ as $n\to\infty$, for all $k\in\mathbb{N}$, implies tightness of $\hat{\mu}(\mathbf{H}_n/\sqrt{n}\,)$ and its weak convergence to the nonrandom measure $\gamma_H$. Since its moments do not depend on the distribution of the i.i.d. sequence $\{X_j\}$, so does the limiting distribution $\gamma_H$, which is symmetric since all its odd moments are zero. By Proposition A.2 it has unbounded support, and is not unimodal. $\square$

4.6. *Markov matrices with centered entries.* In view of Proposition 4.1(ii) we may and shall assume hereafter without loss of generality that the random variables $X_{ij}$ are bounded. Our proof of Theorem 1.3 follows a similar outline as that used in proving Theorems 1.1 and 1.2, where the combinatorial arguments used here rely on matrix decomposition.

Starting with some notation we shall use throughout the proof, let $\Gamma_n$ be a graph whose vertices are two-element subsets of $\{1,\dots,n\}$ with the edges between vertices $a$ and $b$ if the sets overlap, $a\cap b\neq\varnothing$. We indicate that $(a,b)$ is an edge of $\Gamma_n$ by writing $a\sim b$, and for $a\in\Gamma_n$ let $a=\{a^-,a^+\}$ with $1\le a^-<a^+\le n$.

The main tool in the Markov case is the following decomposition:

$$\mathbf{M}_n = \sum_{a\in\Gamma_n} X_a \mathbf{Q}_{a,a},$$

where $X_a := X_{a^+,a^-}$ and $\mathbf{Q}_{a,b}$ is the $n\times n$ matrix defined for vertices $a,b$ of $\Gamma_n$ by

$$\mathbf{Q}_{a,b}[i,j] = \begin{cases} -1, & \text{if } i=a^+, j=b^+, \text{ or } i=a^-, j=b^-,\\ 1, & \text{if } i=a^+, j=b^-, \text{ or } i=a^-, j=b^+,\\ 0, & \text{otherwise.}\end{cases}$$

Let $t_{a,b} = \operatorname{tr}(\mathbf{Q}_{a,b})$. It is straightforward to check that

$$t_{a,b} = \begin{cases} -2, & \text{if } a = b, \\ -1, & \text{if } a \neq b \text{ and } a^- = b^- \text{ or } a^+ = b^+, \\ 1, & \text{if } a^- = b^+ \text{ or } a^+ = b^-, \\ 0, & \text{otherwise.} \end{cases}$$

From this, we see that $t_{a,b} = t_{b,a}$. Since it is easy to check that $\mathbf{Q}_{a,b} \times \mathbf{Q}_{c,d} = t_{b,c}\mathbf{Q}_{a,d}$, we get

$$\operatorname{tr}(\mathbf{Q}_{a_1,a_1} \times \mathbf{Q}_{a_2,a_2} \times \cdots \times \mathbf{Q}_{a_r,a_r}) = \prod_{j=1}^{r} t_{a_j,a_{j+1}}, \tag{4.14}$$

where for convenience we identified $a_{r+1}$ with $a_1$.

For a circuit $\pi = (a_1 \sim \cdots \sim a_r \sim a_1)$ of length $r$ in $\Gamma_n$ let

$$\mathbf{X}_\pi = \prod_{j=1}^{r} t_{a_j,a_{j+1}} \prod_{j=1}^{r} X_{a_j}. \tag{4.15}$$

It follows from (4.14) and (4.15) that

$$\operatorname{tr}(\mathbf{M}_n^r) = \sum_\pi \mathbf{X}_\pi, \tag{4.16}$$

where the sum is over all circuits of length $r$ in $\Gamma_n$, leading to the Markov analog of the path expansion (4.7),

$$\mathbb{E}\operatorname{tr}(\mathbf{M}_n^r) = \sum_\pi \mathbb{E}\mathbf{X}_\pi. \tag{4.17}$$

We say that a circuit $\pi = (a_1 \sim \cdots \sim a_r \sim a_1)$ of length $r$ in $\Gamma_n$ is vertex-matched if for each $i = 1, \ldots, r$ there exists some $j \neq i$ such that $a_i = a_j$, and that it has a match of order 3 if some value is repeated at least three times among $(a_j, j = 1, \ldots, r)$. Note that the only nonvanishing terms in (4.17) come from vertex-matched circuits.

In analogy with Proposition 4.2, we show next that generically vertex-matched circuits have only double repeats, and consequently, the odd moments of $\mathbb{E}\hat{\mu}(\mathbf{M}_n/\sqrt{n})$ converge to zero as $n \to \infty$.

PROPOSITION 4.10. *Fix $r \in \mathbb{N}$. Let $N$ denote the number of vertex-matched circuits in $\Gamma_n$ with $r$ vertices which have at least one match of order 3. Then there is a constant $C_r$ such that for all $n \in \mathbb{N}$*

$$N \leq C_r n^{\lfloor (r+1)/2 \rfloor}.$$

PROOF. Either $r = 2k$ is even, or $r = 2k - 1$ is odd. In both cases, the total number of different vertices per path is at most $k - 1$. Since $a_1 \sim a_2 \sim \cdots \sim a_r$, there are at most $n^2/2$ choices for $a_1$, and then at most $4n$ choices for each of the remaining $k - 2$ distinct values of $a_j$, and one choice for each repeated value. Thus $N \leq 4^r n^2 \times n^{k-2} = Cn^k$. $\square$

Corollary 4.11. *Suppose $\{X_{ij};\, j \geq i \geq 1\}$ are bounded i.i.d. random variables such that $\mathbb{E}(X_{12}) = 0$, $\mathbb{E}(X_{12}^2) = 1$. Then,*

$$\lim_{n\to\infty} \frac{1}{n^{k+1/2}} \mathbb{E}\,\mathrm{tr}(\mathbf{M}_n^{2k-1}) = 0. \tag{4.18}$$

Proof. If $\mathbb{E}\mathbf{X}_\pi$ is nonzero, then all the vertices of the path $a_1 \sim a_2 \sim \cdots \sim a_{2k-1}$ must be repeated at least twice. So for an odd number of vertices, there must be a vertex which is repeated at least three times. Thus, by Proposition 4.10 and the boundedness of $|X_{ij}|$ and of $t_{a,b}$,

$$|\mathbb{E}\,\mathrm{tr}(\mathbf{M}_n^{2k-1})| \leq C_k n^k,$$

and (4.18) follows. □

Let $\mathbf{W}_n = n^{1/2}\mathbf{Z}_n + \mathbf{X}_n + \xi \mathbf{I}_n$, where $\mathbf{X}_n$ is a symmetric $n \times n$ matrix with i.i.d. standard normal random variables (except for the symmetry constraint), $\mathbf{Z}_n = \mathrm{diag}(Z_{ii})_{1\leq i\leq n}$, with i.i.d. standard normal variables $Z_{ii}$ that are independent of $\mathbf{X}_n$ and $\xi$ is a standard normal, independent of all other variables. A direct combinatorial evaluation of the even moments of $\mathbb{E}\hat{\mu}(\mathbf{M}_n/\sqrt{n})$ is provided in [6]. We follow here an alternative, shorter proof, proposed to us by O. Zeitouni. The key step, provided by our next lemma, replaces the even moments by those of the better understood matrix ensemble $\mathbf{W}_n$.

Lemma 4.12. *Suppose $\{X_{ij}; j \geq i \geq 1\}$ is a collection of bounded i.i.d. random variables such that $\mathbb{E}(X_{12}) = 0, \mathbb{E}(X_{12}^2) = 1$. Then, for every $k \in \mathbb{N}$,*

$$\lim_{n\to\infty} n^{-(k+1)}[\mathbb{E}\,\mathrm{tr}(\mathbf{M}_n^{2k}) - \mathbb{E}\,\mathrm{tr}(\mathbf{W}_n^{2k})] = 0. \tag{4.19}$$

Proof. First observe that by Proposition 4.10, we may and shall assume without loss of generality that $\{X_{ij}\}$ is a collection of i.i.d. standard normal random variables, subject to the symmetry constraint $X_{ij} = X_{ji}$ [as such a change affects $n^{-(k+1)}\mathbb{E}\,\mathrm{tr}(\mathbf{M}_n^{2k})$ by at most $C_k n^{-1}$]. Recall the representation $\mathbf{M}_n = \mathbf{X}_n - \mathbf{D}_n$ of (1.3) and let $\widetilde{\mathbf{M}}_n = \mathbf{X}_n - \widetilde{\mathbf{D}}_{n+1}^{(n)}$ where $\widetilde{\mathbf{D}}_{n+1}^{(n)}$ is obtained by omitting the last row and column of the diagonal matrix $\widetilde{\mathbf{D}}_{n+1}$ which is an independent copy of $\mathbf{D}_{n+1}$ that is independent of $\mathbf{X}_n$. Observe that the diagonal entries of $-\widetilde{\mathbf{D}}_{n+1}^{(n)}$ are jointly normal, of zero mean, variance $n+1$ and such that the covariance of each pair is 1. Therefore, with $-\widetilde{\mathbf{D}}_{n+1}^{(n)}$ independent of $\mathbf{X}_n$, for each $n$, the distribution of $\widetilde{\mathbf{M}}_n$ is exactly the same as that of $\mathbf{W}_n$. Consequently, (4.19) is equivalent to

$$\lim_{n\to\infty} n^{-(k+1)}\mathbb{E}[\mathrm{tr}(\mathbf{M}_n^{2k}) - \mathrm{tr}(\widetilde{\mathbf{M}}_n^{2k})] = 0. \tag{4.20}$$

The first step in proving (4.20) is to note that by a path expansion similar to (4.17) we have that

$$\mathbb{E}[\mathrm{tr}(\mathbf{M}_n^{2k}) - \mathrm{tr}(\widetilde{\mathbf{M}}_n^{2k})] = \sum_{\pi} [\mathbb{E}\mathbf{M}_\pi - \mathbb{E}\widetilde{\mathbf{M}}_\pi], \tag{4.21}$$

where now the sum is over all circuits $\pi : \{0, \dots, 2k\} \to \{1, \dots, n\}$, and

$$\mathbf{M}_\pi = \prod_{i=1}^{2k} M_{\pi(i-1),\pi(i)}$$

with the corresponding expression for $\widetilde{\mathbf{M}}_\pi$. Set each word $w$ of length $2k$ to be a circuit by assigning $w[0] = w[2k]$ and let $\Pi(w)$ denote the collection of circuits $\pi$ such that the distinct letters of $w$ are in a one-to-one correspondence with the distinct values of $\pi$. Let $v = v(w)$ be the number of distinct letters in the word $w$, noting that $\#\Pi(w) \le n^{v(w)}$ and that $\mathbb{E}\mathbf{M}_\pi - \mathbb{E}\widetilde{\mathbf{M}}_\pi = f_n(w)$ is independent of the specific choice of $\pi \in \Pi(w)$. Hence, taking the letters of $w$ to be from the set of numbers $\{1, 2, \dots, 2k\}$ with the convention that $w(i) = w[i]$, we identify $w$ as a representative of $\pi \in \Pi(w)$ (recall $w[0] = w[2k]$). For example, $w = abbc$ of $v(w) = 3$ distinct letters becomes $w = 1223$ which we identify with the circuit $\pi \in \Pi(w)$ of length 4 consisting of the edges $\{1, 2\}$, $\{2, 2\}$, $\{2, 3\}$ and $\{3, 1\}$. In view of (4.21), we thus establish (4.20) by showing that for any $w$, some $C_w < \infty$ and all $n$,

$$|f_n(w)| = |\mathbb{E}\mathbf{M}_w - \mathbb{E}\widetilde{\mathbf{M}}_w| \le C_w n^{k-v(w)+1/2}. \tag{4.22}$$

Let $q = q(w)$ be the number of indices $1 \le i \le 2k$ for which $w[i] = w[i-1]$ [e.g., $q(1223) = 1$]. It is clear from the definition of $\mathbf{M}_n$ and $\widetilde{\mathbf{M}}_n$ that $f_n(w) \neq 0$ only if $q(w) \ge 1$. Let $u = u(w)$ count the number of edges of distinct endpoints in $w$, namely, with $\{w[i-1], w[i]\} \in \Gamma_n$, which appear exactly once along the circuit $w$ [e.g., $u(1223) = 3$]. Then, by independence and centering we have that $\mathbb{E}\widetilde{\mathbf{M}}_w = 0$ as soon as $u(w) \ge 1$, whereas it is not hard to check that if $u(w) > q(w)$, then also $\mathbb{E}\mathbf{M}_w = 0$. Thus, it suffices to consider in (4.22) only circuits $w$ with $q(w) \ge u(w)$.

It is not hard to check that excluding the $q$ loop-edges (each connecting some vertex to itself), there are at most $k + \lfloor (u-q)/2 \rfloor$ distinct edges in $w$. These distinct edges form a connected path through $v(w)$ vertices, which for $u \ge 1$ must also be a circuit. Consequently, for any of the words $w$ we are to consider,

$$v(w) \le k + \mathbb{1}_{u(w)=0} + \lfloor (u(w) - q(w))/2 \rfloor \le k. \tag{4.23}$$

Proceeding to bound $|f_n(w)|$, note that any contribution which grows with $n$ must come from the $q$ diagonal entries of $\mathbf{M}_n$ and $\widetilde{\mathbf{M}}_n$ which are

encountered according to the circuit $w$. Suppose first that $u \geq 1$, in which case $f_n(w) = \mathbb{E}\mathbf{M}_w$. Computing the latter, upon expanding the sums in the $q$ relevant diagonal entries of $\mathbf{D}_n = \operatorname{diag}(\sum_{j=1}^n X_{ij})$, we must assign specific choices to at least $u$ of the resulting "free" indices $j_1, \ldots, j_q \in \{1, \ldots, n\}$ in order to match all $u$ unmatched edges of $w$ of the form $\{w[i-1], w[i]\} \in \Gamma_n$. Indeed, by independence and centering, every other term of this expansion has zero expectation. After doing so, as each diagonal entry of $\mathbf{D}_n$ is normal of mean zero and variance $n$, we conclude by Hölder's inequality that $|f_n(w)| \leq C_w n^{(q-u)/2}$. By our bound (4.23) on $v(w)$, this implies that (4.22) holds.

Consider next words $w$ for which $u(w) = 0$ and let $a_1, \ldots, a_q$ be the $q$ vertices for which $\{a_i, a_i\}$ is an edge of the circuit $w$. Let $M_{ii} = Q_i - S_i$ and $\widetilde{M}_{ii} = \widetilde{Q}_i - \widetilde{S}_i$, for $i = 1, \ldots, 2k$, where $Q_i = X_{ii} - \sum_{j=1}^{2k} X_{ij}$, $\widetilde{Q}_i = X_{ii} - \widetilde{X}_{i,n+1} - \sum_{j=1}^{2k} \widetilde{X}_{ij}$ and $\widetilde{S}_i = \sum_{j=2k+1}^n \widetilde{X}_{ij}$ with the corresponding expressions for $S_i$. Note that we may and shall replace each $S_i$ by $\widetilde{S}_i$ without altering $\mathbb{E}\mathbf{M}_w$, and since the off-diagonal entries of $\mathbf{M}_n$ and $\widetilde{\mathbf{M}}_n$ are the same, we have that

$$\begin{aligned} f_n(w) &= \mathbb{E}\left[L_w\left[\prod_{i=1}^q (Q_{a_i} - \widetilde{S}_{a_i}) - \prod_{i=1}^q (\widetilde{Q}_{a_i} - \widetilde{S}_{a_i})\right]\right] \\ &= \sum_{i=1}^q \mathbb{E}\left[L_w (Q_{a_i} - \widetilde{Q}_{a_i}) \prod_{j=1}^{i-1} M_{a_j,a_j} \prod_{j=i+1}^q \widetilde{M}_{a_j,a_j}\right], \end{aligned}$$

where $L_w$ is the product of the $(2k-q)$ off-diagonal entries of $\mathbf{M}_n$ that correspond to the edges of $w$ that are in $\Gamma_n$. Since the distribution of $(L_w, \{Q_i\}, \{\widetilde{Q}_i\})$ is independent of $n > 2k$, while $M_{ii}$ and $\widetilde{M}_{ii}$ are normal of mean zero and variance at most $n+2$, it follows by Hölder's inequality that $|f_n(w)| \leq C_w n^{(q(w)-1)/2}$, which by (4.23) results with (4.22).

As already seen, (4.22) implies that (4.20) holds and hence the proof of the lemma is complete. $\square$

Let $\gamma_0(dx) = \frac{dx}{2\pi}\sqrt{4-x^2}\mathbb{1}_{|x|\leq 2}$ denote the semicircle distribution, $\gamma_1(dx) = \frac{dx}{\sqrt{2\pi}} \exp(-x^2/2)$ denote the standard normal distribution and let $\gamma_M = \gamma_0 \boxplus \gamma_1$ be the corresponding free convolution. In view of Lemma 4.12, our next result shows that the even moments of $\mathbb{E}\hat{\mu}(\mathbf{M}_n/\sqrt{n})$ converge as $n \to \infty$ to those of $\gamma_M$.

PROPOSITION 4.13. *For every $k \in \mathbb{N}$,*

$$\lim_{n\to\infty} n^{-(k+1)} \mathbb{E}\operatorname{tr}(\mathbf{W}_n^{2k}) = \int x^{2k}\, d\gamma_M. \tag{4.24}$$

PROOF. Let $\mathbf{A}_n = \mathbf{Z}_n + n^{-1/2}\xi\mathbf{I}_n$, so $n^{-1/2}\mathbf{W}_n = \mathbf{A}_n + n^{-1/2}\mathbf{X}_n$. By the strong law of large numbers, with probability 1, $\hat{\mu}(\mathbf{A}_n) \to \gamma_1$ weakly. Further, $\sup_n \mathbb{E}\int |x|\,d\hat{\mu}(\mathbf{A}_n) < \infty$, and $\mathbb{E}\int |x|\,d\hat{\mu}(n^{-1/2}\mathbf{X}_n) \le n^{-1}\sqrt{\mathbb{E}\,\mathrm{tr}(\mathbf{X}_n^2)} = 1$, implying by Pastur and Vasilchuk ([17], Theorem 2.1 and page 280), that $\hat{\mu}(\mathbf{W}_n/\sqrt{n})$ converges weakly to $\gamma_M$, in probability. It follows that for any $k \in \mathbb{N}$ and all $r < \infty$,

$$\lim_{n\to\infty} \mathbb{E}\int h_r(x)\,d\hat{\mu}(\mathbf{W}_n/\sqrt{n}) = \int h_r(x)\,d\gamma_M, \tag{4.25}$$

where $h_r(x) = (\min(|x|, r))^{2k}$. Recall that all moments of $\gamma_M$ are finite (cf. Proposition A.3), so as $r \to \infty$ the right-hand side of (4.25) converges to $\int x^{2k}\,d\gamma_M$. It is not hard to check that for any $k \in \mathbb{N}$,

$$\mathbb{E}\int x^{2k}\,d\hat{\mu}(\mathbf{W}_n/\sqrt{n}) = n^{-(k+1)}\mathbb{E}\,\mathrm{tr}(\mathbf{W}_n^{2k})$$

is bounded in $n$ by some $C_k < \infty$. Hence, for all $n$,

$$\left| n^{-(k+1)}\mathbb{E}\,\mathrm{tr}(\mathbf{W}_n^{2k}) - \mathbb{E}\int h_r(x)\,d\hat{\mu}(\mathbf{W}_n/\sqrt{n}) \right| \le C_{k+1} r^{-2},$$

and (4.24) follows by considering $r \to \infty$ in (4.25). $\square$

We next derive the analog of Proposition 4.3 and similarly to Proposition 4.9, get as a result the concentration of moments of $\hat{\mu}(\mathbf{M}_n/\sqrt{n})$ around those of $\mathbb{E}(\hat{\mu}(\mathbf{M}_n/\sqrt{n}))$.

PROPOSITION 4.14. *Fix $r \in \mathbb{N}$. Let $N$ denote the number of vertex-matched quadruples of circuits in $\Gamma_n$ with $r$ vertices each, such that none of them is self-matched. Then there is a constant $C_r$ such that*

$$N \le C_r n^{2r+2}.$$

PROOF. Let $\mathcal{P}$ denote the partition of the $4r$ vertices of the circuits $\pi_1, \dots, \pi_4$ in $\Gamma_n$ to $|\mathcal{P}| \le 2r$ distinct groups of matching vertices, with at least two elements in each group, while having each circuit cross-matched to at least one of the other circuits. As part of $\mathcal{P}$ we specify also which of the four types of edges to use in each connection along the circuits. For $i = 1, 2, 3, 4$, let $u_i = u_i(\mathcal{P})$ be the number of distinct vertices in $\pi_i$ that do not appear in any $\pi_j$, $j < i$. There are at most $n^{1+u_1}$ ways to choose the circuit $\pi_1$ in agreement with $\mathcal{P}$, that is, $n^2/2$ ways to choose the vertex $a_1$ of $\pi_1$ and at most $n$ ways for each of the remaining $u_1 - 1$ distinct vertices of $\pi_1$. For $i = 2, 3, 4$, per given $\pi_j$, $j < i$, the same procedure shows that there are at most $n^{1+u_i}$ ways to complete the circuit $\pi_i$. Further, if $\pi_i$ is cross-matched to $\pi_j$ for some $j < i$, then starting the completion of $\pi_i$ at a vertex that we already determined by such a cross-match, we have that there are only $n^{u_i}$ ways to complete $\pi_i$. The latter improved bound always applies for $i=4$, and it is not hard to check that upon re-ordering the four circuits, we can assure that it applies also for $i=3$. We thus get at most $n^{u+2}$ quadruples of circuits per choice of $\mathcal{P}$, where $u=\sum_i u_i=|\mathcal{P}|\le 2r$, yielding the stated bound. □

PROPOSITION 4.15. *Suppose $\{X_{ij};\, j\ge i\ge 1\}$ is a collection of bounded i.i.d. random variables such that $\mathbb{E}(X_{12})=0$ and $\mathbb{E}(X_{12}^2)=1$. For any $r\in\mathbb{N}$, there exists $C_r<\infty$ such that $\mathbb{E}[(\operatorname{tr}(\mathbf{M}_n^r)-\mathbb{E}\operatorname{tr}(\mathbf{M}_n^r))^4]\le C_r n^{2r+2}$ for all $n\in\mathbb{N}$.*

PROOF. By (4.16) we have the Markov analog of (4.13)

$$\mathbb{E}[(\operatorname{tr}(\mathbf{M}_n^r)-\mathbb{E}\operatorname{tr}(\mathbf{M}_n^r))^4]=\sum_{\pi_1,\pi_2,\pi_3,\pi_4}\mathbb{E}\left[\prod_{j=1}^4(\mathbf{X}_{\pi_j}-\mathbb{E}(\mathbf{X}_{\pi_j}))\right], \tag{4.26}$$

where the sum is taken over all circuits $\pi_j$, $j=1,\dots,4$, in $\Gamma_n$, each having $r$ vertices. With the random variables $\{X_{ij};n\ge j\ge i\ge 1\}$ independent and of mean zero, just like the proof of Proposition 4.9, it suffices to take the sum in (4.26) over all vertex-matched quadruples of circuits on $\Gamma_n$, such that none of them is self-matched. Since $|X|$ (and hence $|\mathbf{X}_\pi|$) is bounded the stated inequality follows from the bound of Proposition 4.14 on the number of such quadruples. □

PROOF OF THEOREM 1.3. The proof is very similar to that of Theorems 1.1 and 1.2, where by Proposition 4.1(ii), we may and shall assume that $\{X_{ij};j\ge i\ge 1\}$ is a collection of i.i.d. bounded random variables. Then, by (4.18) the odd moments of the average measure $\mathbb{E}(\hat\mu(\mathbf{M}_n/\sqrt{n}))$ converge to 0, and by Proposition 4.13 the even moments converge to those of $\gamma_M$, whereas from Proposition 4.15 we conclude that with probability 1 the same applies to the moments of $\hat\mu(\mathbf{M}_n/\sqrt{n})$. By Proposition A.3, $\gamma_M$ is a symmetric measure of bounded smooth density that, though of unbounded support, is uniquely determined by its moments (having in particular zero odd moments). Hence, the almost sure convergence $\int x^k\,d\hat\mu(\mathbf{M}_n/\sqrt{n})\to\int x^k\,d\gamma_M$ as $n\to\infty$, for all $k\in\mathbb{N}$, implies the weak convergence of $\hat\mu(\mathbf{M}_n/\sqrt{n})$ to $\gamma_M$. □

## APPENDIX

**A.1. Properties of $\gamma_H$, $\gamma_M$ and $\gamma_T$.** In this section we establish properties of the symmetric measures with moments given by (3.5), and (3.8) and the free convolution $\gamma_M$ of Theorem 1.3. For proofs, it is convenient to express the volumes $p_H(w)$ and $p_T(w)$ as the probabilities that involve sums of independent uniform random variables. This can be done by setting the undetermined variables as the independent uniform $U[0,1]$ random variables $U_0, U_1, \ldots, U_k$, expressing the dependent variables as the linear combinations of $U_0, U_1, \ldots, U_k$, and expressing the volumes as the probabilities that these linear combinations are in the interval $I$. For each partition word $w$ of length $2k$ with a nonzero volume $p(w)$, this probability takes the form

$$p(w) = \mathbb{P}\left( \bigcap_{i=1}^{k} \left\{ \sum_{j=0}^{M} n_{i,j} U_j \in [0,1] \right\} \right), \tag{A.1}$$

where $n_{i,j}$ are integers and $M = k$.

Proposition A.1. *A symmetric measure $\gamma_T$ with even moments given by* (3.5) *has unbounded support.*

Proof. It suffices to show that $(m_{2k})^{1/k} \to \infty$. Let $w$ be a partition word of length $2k$. Denoting $S_i = \sum_j n_{i,j} U_j - \frac{1}{2}, i = 1, 2, \ldots, k$, we have

$$p_T(w) = \mathbb{P}\left( \bigcap_{i=1}^{k} \{ |S_i| < \tfrac{1}{2} \} \right). \tag{A.2}$$

Since the coefficients $n_{i,j}$ in (A.1) take values $0, \pm 1$ only, and $\sum_j n_{i,j} = 1$, each of the sums $S_i$ in (A.2) has the following form:

$$S = (U_\alpha - 1/2) + \sum_{j=1}^{L} (U_{\beta(j)} - U_{\gamma(j)}), \tag{A.3}$$

where $\alpha, \beta(j), \gamma(j), j = 1, \ldots, L$, are all different. Let $L_i$ denote the number of independent random variables $U$ in this representation for $S_i$. Clearly, $1 \le L_i \le k+1$.

Fixing $\varepsilon > 0$ let $U_j = 1/2 + V_j/(\varepsilon(k+1))$ for $j = 0, \ldots, k$. For $k > 1/\varepsilon$ define the event

$$A = \bigcap_{j=0}^{k} \left\{ |U_j - 1/2| < \frac{1}{2\varepsilon(k+1)} \right\},$$

noting that conditionally on $A$, the random variables $V_0, \ldots, V_k$ are independent, each uniformly distributed on $[-1/2, 1/2]$. As under this conditioning the i.i.d. random variables $\{V_j\}$ have symmetric laws, it is easy to check that for $i = 1, \ldots, k$, the form (A.3) of $S_i$ implies that

$$\mathbb{P}(|S_i| > \tfrac{1}{2} | A) = \mathbb{P}\left( \left| \sum_{j=1}^{L_i} V_j \right| > \varepsilon(k+1)/2 \right) = 2\mathbb{P}\left( \sum_{j=1}^{L_i} V_j > \varepsilon(k+1)/2 \right),$$

which by Markov's inequality is bounded above by

$$2e^{-\varepsilon^2(k+1)/2}(\mathbb{E}e^{\varepsilon V})^{L_i} = e^{-\varepsilon^2(k+1)/2}\left(\frac{e^{\varepsilon/2} - e^{-\varepsilon/2}}{\varepsilon}\right)^{L_i}.$$

Since $\frac{e^x - e^{-x}}{2x} \le e^{x^2/2}$ for $x > 0$, and $L_i \le k+1$, we deduce that

$$\mathbb{P}(|S_i| > \tfrac{1}{2}|A) \le 2\exp(-\varepsilon^2(k+1)/2 + \varepsilon^2 L_i/4) \le 2e^{-\varepsilon^2(k+1)/4}, \tag{A.4}$$

for $i = 1, \dots, k$. As $2ke^{-\varepsilon^2(k+1)/4} \le 1/2$ for some $k_0 = k_0(\varepsilon) < \infty$ and all $k \ge k_0$, it follows from (A.2) and (A.4) that for all $k \ge k_0$ and any word $w$ of length $2k$,

$$p_T(w) \ge \tfrac{1}{2}\mathbb{P}(A) = \tfrac{1}{2}(\varepsilon(k+1))^{-(k+1)}. \tag{A.5}$$

Since there are more than $k!$ partition words $w$ of length $2k$, this shows that for all large enough $k$ we have

$$m_{2k} \ge \tfrac{1}{2}k!(\varepsilon(k+1))^{-(k+1)} \ge (3\varepsilon)^{-k}.$$

Hence, $\limsup_{k\to\infty} m_{2k}^{1/k} \ge 1/(3\varepsilon)$. As $\varepsilon > 0$ is arbitrarily small, this completes the proof. □

PROPOSITION A.2. *A symmetric measure $\gamma_H$ with even moments given by* (3.8) *is not unimodal and has unbounded support.*

PROOF. Suppose that the symmetric distribution $\gamma_H$ is unimodal. Since all moments of $\gamma_H$ are finite, from Khinchin's theorem (see [14], Theorem 4.5.1), it follows that if $\phi(t) = \int e^{itx}\gamma_H(dx)$ denotes the characteristic function of $\gamma_H$, then $g(t) = \phi(t) + t\phi'(t)$ must be a characteristic function, too. The even moments corresponding to $g(t)$ are $(2k+1)m_{2k}(\gamma_H)$, and must be a positive definite sequence, that is, the Hankel matrices with entries $[(2(i+j)-3)m_{2(i+j-2)}(\gamma_H)]_{1\le i,j\le n}$ should all be nonnegative definite. However, with $m_4 = 2$, $m_6 = 11/2$ and $m_8 = 281/15$, for $n = 3$ the determinant

$$\det\begin{bmatrix} 1 & 3m_2 & 5m_4 \\ 3m_2 & 5m_4 & 7m_6 \\ 5m_4 & 7m_6 & 9m_8 \end{bmatrix} = \det\begin{bmatrix} 1 & 3 & 10 \\ 3 & 10 & 77/2 \\ 10 & 77/2 & 843/5 \end{bmatrix} = -73/20$$

is negative. Thus, $\gamma_H$ is not unimodal.

To show that the support of $\gamma_H$ is unbounded we proceed like in the Toeplitz case. The main technical obstacle is that some partition words contribute zero volume. We will therefore have to find enough partition words that contribute a nonzero volume, and then give a lower bound for this contribution.

We consider only moments of order $4k-2$, $k \geq 2$, and find the contribution of the partition words which have no repeated letters in the first half, that is,

$$w[1] \neq w[2] \neq \cdots \neq w[2k-1].$$

That is, we consider the set of partition words $w$ of length $4k-2$ of the form $w = abc\ldots$ with the first $2k-1$ letters written in the fixed (alphabetic) order, followed by the repeated letters $a, b, c, \ldots$ at positions $2k, \ldots, 4k-2$. We also require that the repeats are placed at odd distance from the original matching letter. Formally, we consider the set of partition words $w$ of length $4k-2$ which satisfy the following condition:

$$\text{If } w[\alpha] = w[\beta] \text{ and } \alpha < \beta, \qquad \text{then } \alpha \not\equiv \beta \bmod 2,\ \alpha \leq 2k-1 \text{ and } \beta \geq 2k.$$

Since we can permute all letters at locations $2k, 2k+2, \ldots, 4k-2$ and all letters at locations $2k+1, 2k+3, \ldots, 4k-3$, clearly there are $k!(k-1)!$ such partition words.

To show that all such partition words contribute a nonzero volume, we need to carefully analyze the matrix of the resulting system of equations (3.6). This is a $(2k-1) \times (4k-1)$ matrix with entries $0, \pm 1$ only. The first $2k-1$ columns of the matrix are filled in with the pattern of sliding pairs $1, 1$ corresponding to first occurrences of every letter, that is, the left-hand sides of equations (3.6) are simply

$$\begin{cases} x_0 + x_1 & = \ldots \\ \quad x_1 + x_2 & = \ldots \\ & \vdots \\ \qquad x_{2k-2} + x_{2k-1} & = \ldots. \end{cases}$$

So the first $2k$ columns of the matrix are as follows, with the star denoting as yet unspecified entries of the $2k$th column.

```
1100..00*
0110..00*
0011..00*
 ...
0000..11*
0000..011
```

The remaining columns are as follows. In every even row of the second half we have disjoint (nonoverlapping) pairs $(-1, -1)$, including the site adjacent to the "last letter," that has entry 1 in the last row, and entry $-1$ in one of the odd rows. None of these $-1, -1$ are in the last column, a coefficient of $x_{4k-2}$.

In the odd rows we have pairs of consecutive $(-1,-1)$ which overlap entries from the even rows, but not themselves, including a single $(-1,-1)$ pair which fills in one spot in the last column, the coefficients of $x_{4k-2}$.

For example, the word $w = abc\dots abc\dots$, where all $2k-1$ letters $a, b, c, \dots$ are repeated alphabetically twice, is in the class of the partition words under consideration. The corresponding system of equations is

$$\begin{cases} x_0 + x_1 = x_{2k-1} + x_{2k} & \\ \vdots & \\ x_i + x_{i+1} = x_{2k+i-1} + x_{2k+i}, & i = 1, 2, \dots, 2k-3, \\ \vdots & \\ x_{2k-2} + x_{2k-1} = x_{4k-3} + x_{4k-2}, & \end{cases}$$

and its matrix is

```
1100..00-1-1 0 ... 0 0
0110..00 0-1-1 ... 0 0
0011..00 0 0-1 ... 0 0
...
0000..11 0 0 0 ...-1 0
0000..01 1 0 0 ...-1-1
```

All other partition words in our class are obtained from permuting letters $w[2k], w[2k+2], \dots, w[4k-2]$, and then permuting letters $w[2k+1], w[2k+3], \dots, w[4k-3]$ of $w = abc\dots abc$. Thus all other systems of equations are obtained from the above one by permuting even rows in columns $2k+1, 2k+2, \dots, 4k-2$ and odd rows in columns $2k, 2k+1, \dots, 4k-1$ (apart from the 1 at column $2k$ and row $2k-1$ which is never permuted, but gets eliminated if the first row permutes to become the last one). For each of these words the sum of all odd rows in the system minus the sum of all even rows is $[1, 0, \dots, 0, -1]$, implying that for such $w$ the additional constraint $x_0 = x_{4k-2}$ we require when computing $p_H(w)$ is merely a consequence of (3.6).

The solutions of equations (3.6) for such partition words $w$ are easy to analyze due to parity considerations. Gaussian elimination consists here of subtractions of the given row from the row directly above it, starting with the subtraction of the $(2k-1)$ row and ending with the subtraction of the second row from the first row, at which point the first $2k-1$ columns become the identity matrix. During these subtractions, a $-1$ entry in each column of the original system can meet a nonzero entry only from a row positioned at an odd distance above it, in which case they cancel each other. So as we keep subtracting, all coefficients take values $0, \pm 1$ only. Further, for each row the sum of the entries in columns $2k, \dots, 4k-1$ is $-2$, except for the last row for which it is $-1$. Thus, after all subtractions have been made, these sums are $-1$ at each of the rows. We can now set the $2k$ undetermined variables to i.i.d. $U[0,1]$ random variables, $x_{2k-1} = U_0, \ldots, x_{4k-2} = U_{2k-1}$, and solve the $2k-1$ equations for the dependent variables $x_0, \ldots, x_{2k-2}$. By the above considerations we know that each of these dependent random variables is expressed as an alternating sum of independent uniform $U[0,1]$ random variables of the form (A.3).

The argument we used for deriving (A.5) thus gives the bound $p_H(w) \geq \frac{1}{2}(2k\varepsilon)^{-2k}$ for each of these $k!(k-1)!$ partition words, and hence for all $k$ large enough, we have

$$m_{4k-2}(\gamma_H) \geq \tfrac{1}{2}k!(k-1)!(2\varepsilon k)^{-2k} \geq (6\varepsilon)^{-2k}.$$

Thus $m_{4k-2}^{1/k} \to \infty$, which implies that the support of $\gamma_H$ is unbounded. □

PROPOSITION A.3. *The free convolution $\gamma_M = \gamma_0 \boxplus \gamma_1$ of the standard semicircle distribution $\gamma_0$ and the standard normal $\gamma_1$ is a symmetric measure, determined by moments, has unbounded support and a smooth bounded density.*

PROOF. By Corollary 2 in [2], $\gamma_M$ has a density, by Corollary 4 in [2] the density is smooth and by Proposition 5 in [2] it is bounded.

We now verify that $\gamma_M$ is determined by moments and has unbounded support. We need the following observation: *a probability measure $\mu$ has odd moments vanishing iff the odd free cumulants $k_{2r+1}(\mu)$ of $\mu$ vanish.* This can be easily read from formula (72) in [22].

Since free cumulants linearize the free convolution, $k_r(\gamma_M) = k_r(\gamma_0) + k_r(\gamma_1)$. This shows that the odd moments of $\gamma_M$ vanish. Recall that the free cumulants $k_n(\mu)$ and the moments $m_n(\mu)$ of a probability measure $\mu$ are related by formula (72) in [22]. In particular, for $\mu$ with vanishing odd moments, the even cumulants $k_{2r}(\mu)$ are related to the moments by the equations

$$m_{2n}(\mu) = \sum_{r=1}^{n} k_{2r}(\mu) \sum_{i_1+\cdots+i_{2r}=2n-2r} \prod_{j=1}^{2r} m_{i_j}(\mu), \qquad n = 1, 2, \ldots. \tag{A.6}$$

By symmetry, the odd cumulants of $\gamma_1$ vanish, and $k_{2r}(\gamma_1)$ are nonnegative; $k_{2r}(\gamma_1)$ count all irreducible pair partitions of $\{1, \ldots, 2r\}$ (see [5], page 152). Since $k_2(\gamma_0) = 1$, and all higher free cumulants of $\gamma_0$ vanish (see [12], Example 2.4.6), we have

$$k_{2r}(\gamma_1) \leq k_{2r}(\gamma_M) \leq 2k_{2r}(\gamma_1).$$

Together with (A.6) this implies by induction that

$$m_{2r}(\gamma_1) \leq m_{2r}(\gamma_M) \leq 4^r m_{2r}(\gamma_1).$$

In particular, $\gamma_M$ has unbounded support and is uniquely determined by moments. Since its odd cumulants vanish, the odd moments vanish and $\gamma_M$ is symmetric. □

**A.2. Moments of free convolution.** In this section we identify moments of the free convolution $\gamma_0 \boxplus \gamma_1$. The result and the method of proof were suggested by Bożejko and Speicher [5], who give a combinatorial expression for the moments of free convolutions of normal densities.

Denote by $\mathcal{W}$ the set of all partition words. Recall that a (partition) subword of a word $w$ is a partition word $w_1$ such that $w = a\ldots cw_1d..z$. Let $\mathcal{W}_0$ be the set of all *irreducible partition words*, that is, words that have no proper (nonempty) partition subwords.

Definition A.1 ([5]). We say that $p:\mathcal{W}\to\mathbb{R}$ is *pyramidally multiplicative*, if for every $w\in\mathcal{W}$ of the form $w = a\ldots cw_1d..z$, we have $p(w) = p(w_1)p(a\ldots cd..z)$.

Lemma A.4 ([5], page 152). *Suppose that the moments are given by*

$$m_{2n} = \sum_{w\in\mathcal{W},|w|=2n} p(w), \tag{A.7}$$

*and* $m_{2n-1} = 0$, $n = 1, 2, \ldots.$ *If the weights* $p(w)$ *are pyramidally multiplicative, then the free cumulants are*

$$k_{2n} = \sum_{w\in\mathcal{W}_0,|w|=2n} p(w).$$

Proposition A.5. *A symmetric measure* $\gamma_M$ *with the even moments given by* (3.1) *is given by the free convolution* $\gamma_M = \gamma_0 \boxplus \gamma_1$.

Proof. We apply Lemma A.4 to measures $\gamma_M$, $\gamma_0$ and $\gamma_1$. If $w = ..w_1..$, then $h(w) = h(w_1) + h(w\setminus w_1)$, so the Markov weights $p_M(w) := 2^{h(w)}$ are pyramidally multiplicative. It is well known that the moments of the normal distribution are given by (A.7) with $p_1(w) = 1$, which is (trivially) multiplicative. The moments of the semicircle distribution are given by (A.7) with $p_0(w) = 1$ for the so-called noncrossing words, and $p_0(w) = 0$ otherwise. (A partition word is noncrossing, if it can be reduced to the empty word by removing pairs of consecutive double letters $xx$, one at a time.) It is well known that this weight is pyramidally multiplicative, too.

We now use Lemma A.4 to compare the free cumulants of the semicircle, normal and Markov distributions. Let $w\in\mathcal{W}_0$. If $|w| = 2$, then $p_M(w) = 2$, and otherwise $p_M(w) = 2^0 = 1$ as an irreducible word has no proper subwords, and hence no encapsulated subwords. Thus $k_2(\gamma_M) = 2$, and for $n\ge 2$

$$k_{2n}(\gamma_M) = \#\{w\in\mathcal{W}_0, |w| = 2n\}.$$

If $|w| = 2$, then $p_0(aa) = 1$, and otherwise $p_0(w) = 0$ as an irreducible word of length 4 or more cannot be noncrossing. Thus $k_2(\gamma_0) = 1$, and for $n\ge 2$

$$k_{2n}(\gamma_0) = 0.$$

From $p_1(w) = 1$ we get

$$k_{2n}(\gamma_1) = \#\{w \in \mathcal{W}_0, |w| = 2n\}$$

for $n \geq 1$; in particular, $k_2(\gamma_1) = 1$. Thus, for $n \geq 1$

$$k_{2n}(\gamma_M) = k_{2n}(\gamma_0) + k_{2n}(\gamma_1),$$

which proves that $\gamma_M = \gamma_0 \boxplus \gamma_1$. $\square$

**Acknowledgments.** Part of the research of WB was conducted while visiting the Department of Statistics of Stanford University. The authors thank Marek Bożejko, Persi Diaconis, J. T. King, Rafał Latała, Qiman Shao, Ronald Speicher and Richard P. Stanley for helpful comments, references and encouragement, Ofer Zeitouni for a shorter proof of Theorem 1.3 and additional comments, A. Sakhanenko for electronic access to his papers, and Steven Miller, Chris Hammond and Arup Bose for information about their research.

## REFERENCES

[1] Bai, Z. D. (1999). Methodologies in spectral analysis of large-dimensional random matrices, a review. *Statist. Sinica* **9** 611–677. MR1711663

[2] Biane, P. (1997). On the free convolution with a semi-circular distribution. *Indiana Univ. Math. J.* **46** 705–718. MR1488333

[3] Bose, A., Chatterjee, S. and Gangopadhyay, S. (2003). Limiting spectral distributions of large dimensional random matrices. *J. Indian Statist. Assoc.* **41** 221–259. MR2101995

[4] Bose, A. and Mitra, J. (2002). Limiting spectral distribution of a special circulant. *Statist. Probab. Lett.* **60** 111–120. MR1945684

[5] Bożejko, M. and Speicher, R. (1996). Interpolations between bosonic and fermionic relations given by generalized Brownian motions. *Math. Z.* **222** 135–159. MR1388006

[6] Bryc, W., Dembo, A. and Jiang, T. (2003). Spectral measure of large random Hankel, Markov and Toeplitz matrices. Expanded version available at http://arxiv.org/abs/math.PR/0307330.

[7] Diaconis, P. (2003). Patterns in eigenvalues: The 70th Josiah Willard Gibbs lecture. *Bull. Amer. Math. Soc. (N.S.)* **40** 155–178 (electronic). MR1962294

[8] Dudley, R. M. (2002). *Real Analysis and Probability*. Cambridge Univ. Press. MR1932358

[9] Fulton, W. (2000). Eigenvalues, invariant factors, highest weights, and Schubert calculus. *Bull. Amer. Math. Soc. (N.S.)* **37** 209–249 (electronic). MR1754641

[10] Grenander, U. and Szegő, G. (1984). *Toeplitz Forms and Their Applications*, 2nd ed. Chelsea, New York. MR890515

[11] Hammond, C. and Miller, S. (2005). Eigenvalue density distribution for real symmetric Toeplitz ensembles. *J. Theoret. Probab.* **18** 537–566.

[12] Hiai, F. and Petz, D. (2000). *The Semicircle Law, Free Random Variables and Entropy*. Amer. Math. Soc., Providence, RI. MR1746976

[13] Lidskiĭ, V. B. (1950). On the characteristic numbers of the sum and product of symmetric matrices. *Dokl. Akad. Nauk SSSR (N.S.)* **75** 769–772. MR39686

[14] Lukacs, E. (1970). *Characteristic Functions*, 2nd ed. Hafner, New York. MR346874

[15] Mohar, B. (1991). The Laplacian spectrum of graphs. In *Graph Theory, Combinatorics, and Applications* **2** 871–898. Wiley, New York. MR1170831

[16] Nicolas, J.-L. (1992). An integral representation for Eulerian numbers. In *Sets, Graphs and Numbers. Colloq. Math. Soc. János Bolyai* **60** 513–527. North-Holland, Amsterdam. MR1218216

[17] Pastur, L. and Vasilchuk, V. (2000). On the law of addition of random matrices. *Comm. Math. Phys.* **214** 249–286. MR1796022

[18] Sakhanenko, A. I. (1985). Estimates in an invariance principle. In *Limit Theorems of Probability Theory. Trudi Inst. Math.* **5** 27–44, 175. Nauka, Novosibirsk. MR821751

[19] Sakhanenko, A. I. (1991). On the accuracy of normal approximation in the invariance principle. *Siberian Adv. Math.* **1** 58–91. MR1138005

[20] Sen, A. and Srivastava, M. (1990). *Regression Analysis*. Springer, New York. MR1063855

[21] Serre, J.-P. (1997). Répartition asymptotique des valeurs propres de l'opérateur de Hecke $T_p$. *J. Amer. Math. Soc.* **10** 75–102. MR1396897

[22] Speicher, R. (1997). Free probability theory and non-crossing partitions. *Sém. Lothar. Combin.* **39** Art. B39c (electronic). MR1490288

[23] Tanny, S. (1973). A probabilistic interpretation of Eulerian numbers. *Duke Math. J.* **40** 717–722. MR340045

[24] Wigner, E. P. (1958). On the distribution of the roots of certain symmetric matrices. *Ann. of Math. (2)* **67** 325–327. MR95527

W. Bryc
Department of Mathematics
University of Cincinnati
P.O. Box 210025
Cincinnati, Ohio 45221
USA
e-mail: wlodzimierz.bryc@uc.edu
url: math.uc.edu/~brycw/

A. Dembo
Department of Statistics
and Department of Mathematics
Stanford University
Stanford, California 94305
USA
e-mail: amir@math.stanford.edu
url: www-stat.stanford.edu/~amir/

T. Jiang
School of Statistics
313 Ford Hall
224 Church Street S.E.
Minneapolis, Minnesota 55455
USA
e-mail: tjiang@stat.umn.edu
url: www.stat.umn.edu/~tjiang/